\documentclass[oneside,english,12pt]{amsart}
\usepackage{lmodern}
\usepackage[T1]{fontenc}       
\usepackage[latin9]{inputenc} 
\usepackage[english]{babel}   
\usepackage[top=3cm, bottom=3cm, left=3cm, right=3cm, heightrounded, marginparwidth=2.5cm, marginparsep=6mm]{geometry}

\usepackage{setspace} 
\usepackage{marginnote}
\usepackage{verbatim}
\usepackage{mathrsfs}	
\usepackage{amstext}
\usepackage{amsthm}
\usepackage{amssymb}
\usepackage{amsopn} 
\usepackage{bbold}
\usepackage{stix}
\usepackage{enumitem}				
\usepackage[all]{xy}    				
\usepackage{mathdots}   				
\usepackage{aliascnt}   				
\usepackage{esint}					
\usepackage{tikz}  \usetikzlibrary{matrix}	
\usepackage[backref=page]{hyperref}
\usepackage{mathdots}
\hypersetup{
    colorlinks,
    linkcolor={red!50!black},
    citecolor={green!50!black},
    urlcolor={blue!80!black},
    linktocpage
}
\usepackage{soul}  
\usepackage{longtable} 
\usepackage{lipsum}
\usepackage{footnote}
\makesavenoteenv{tabular}

\usepackage{multirow}
\usepackage{longtable} 				
\allowdisplaybreaks

\renewcommand*{\backref}[1]{}
\renewcommand*{\backrefalt}[4]{%
    \ifcase #1 (Not cited.)%
    \or        (Cited on page~#2.)%
    \else      (Cited on pages~#2.)%
    \fi}

\numberwithin{equation}{section}
\numberwithin{figure}{section}

\theoremstyle{theorem}
  \newtheorem*{cor*}{Corollary}
  \newtheorem*{conj*}{Isomorphism Conjecture}
  \newtheorem*{thm*}{Theorem}
  \newtheorem*{lem*}{Lemma}
  \newtheorem*{claim*}{Claim}
  \newtheorem*{mclaim*}{Main Claim}
  
  \newtheorem{thmx}{Theorem}

  \newaliascnt{corx}{thmx}
  \newtheorem{corx}[corx]{Corollary}
  \aliascntresetthe{corx}

  \newtheorem{thm}{Theorem}[section]

  \newaliascnt{lem}{thm}
  \newtheorem{lem}[lem]{Lemma}
  \aliascntresetthe{lem}

  \newaliascnt{klem}{thm}
  
  \aliascntresetthe{klem}

  \newaliascnt{cor}{thm}
  \newtheorem{cor}[cor]{Corollary}
  \aliascntresetthe{cor}

  \newaliascnt{prop}{thm}  
  \newtheorem{prop}[prop]{Proposition}
  \aliascntresetthe{prop}

\theoremstyle{definition}
  \newaliascnt{defn}{thm}
  \newtheorem{defn}[defn]{Definition}
  \aliascntresetthe{defn}

  \newaliascnt{exmpl}{thm}
  \newtheorem{exmpl}[exmpl]{Example}
  \aliascntresetthe{exmpl}

\theoremstyle{remark}
  \newaliascnt{rem}{thm}
  \newtheorem{rem}[rem]{Remark}
  \aliascntresetthe{rem}

  \theoremstyle{remark}
  \newaliascnt{con}{thm}
  \newtheorem{con}[con]{Construction}
  \aliascntresetthe{con}

\newcommand{\bbC}{\mathbb{C}}

\newcommand{\bbH}{\mathbb{H}}

\newcommand{\bbK}{\mathbb{K}}

\newcommand{\bbN}{\mathbb{N}}

\newcommand{\bbP}{\mathbb{P}}

\newcommand{\bbR}{\mathbb{R}}             
\newcommand{\bbS}{\mathbb{S}}

\newcommand{\bbZ}{\mathbb{Z}}      
\newcommand{\bbone}{\mathbb{1}} 

\newcommand{\clD}{\mathcal{D}}
\newcommand{\clE}{\mathcal{E}}
\newcommand{\clF}{\mathcal{F}}
\newcommand{\clG}{\mathcal{G}}
\newcommand{\clH}{\mathcal{H}}

\newcommand{\clJ}{\mathcal{J}}

\newcommand{\clO}{\mathcal{O}}
\newcommand{\clP}{\mathcal{P}}

\newcommand{\clT}{\mathcal{T}}

\newcommand{\clV}{\mathcal{V}}

\newcommand{\scC}{\mathscr{C}}
\newcommand{\scD}{\mathscr{D}}

\newcommand{\scL}{\mathscr{L}}

\newcommand{\scS}{\mathscr{S}}

\newcommand{\kay}{\mathscr{k}}

\DeclareMathOperator{\Aut}{Aut}

\DeclareMathOperator{\Char}{char}
\DeclareMathOperator{\dd}{d}

\DeclareMathOperator{\ev}{ev}
\DeclareMathOperator{\Gr}{Gr}    
\DeclareMathOperator{\id}{id}          
\DeclareMathOperator{\HH}{H}
\DeclareMathOperator{\Hom}{Hom}
\DeclareMathOperator{\im}{im}

\DeclareMathOperator{\Prob}{Prob}

\DeclareMathOperator{\sgn}{sgn}
\DeclareMathOperator{\Span}{span} 

\DeclareMathOperator{\St}{st}

\DeclareMathOperator{\Sym}{Sym}

\newcommand{\Hb}{{\rm H}_{\rm b}}

\newcommand{\Linfty}{L^\infty}

\newcommand{\bcdot}{{\scriptscriptstyle \bullet}}

\newcommand{\BX}{\bar{X}}

\newcommand{\C}{\mathbb C}
\newcommand{\mres}{\, \mathbin{\vrule height 1.6ex depth 0pt width
0.13ex\vrule height 0.13ex depth 0pt width 1.3ex}}

\DeclareMathOperator{\GL}{GL}         
\DeclareMathOperator{\SL}{SL}         
\DeclareMathOperator{\OO}{O}               
\DeclareMathOperator{\SO}{SO}             
\DeclareMathOperator{\UU}{U}               
\DeclareMathOperator{\SU}{SU}             
\DeclareMathOperator{\Sp}{Sp}              

\newcommand{\qand}{\quad \mathrm{and} \quad}
\newcommand{\qqand}{\qquad \mathrm{and} \qquad}

\AtBeginDocument{\addtocontents{toc}{\protect\setlength{\parskip}{0pt}}}

\begin{document}

\title[Stabilization of bounded cohomology]{Stabilization of bounded cohomology for classical groups}
\author{Carlos De la Cruz Mengual}
\address{Faculty of Electrical and Computer Engineering \\ Technion, Haifa, Israel}
\email{c.delacruz@technion.ac.il}

\author{Tobias Hartnick}
\address{Institute of Algebra and Geometry \\
KIT, Karlsruhe, Germany}
\email{tobias.hartnick@kit.edu}

\begin{abstract}
We show that bounded cohomology stabilizes along sequences of classical Lie groups, and along sequences of lattices in them. Our method is based on a criterion from \cite{DMH0} which adapts Quillen's stability method to the setting of bounded cohomology. This criterion is then applied to a family of measured complexes, the so-called Stiefel complexes, associated to any vector space endowed with a non-degenerate sesquilinear form.
\end{abstract}

\maketitle
\vspace{-15pt}

\section{Introduction}
\subsection{Bounded-cohomological stability and the isomorphism conjecture}
The present article is part of a series of articles \cite{DMH0, DM} which is concerned with the explicit computation of the continuous bounded cohomology of Lie groups using stability results. \emph{Continuous bounded cohomology} (as introduced by Burger and Monod \cite{Burger-Monod1, Burger-Monod2, Monod-Book} in the early 2000s) is a generalization of ordinary bounded cohomology (as introduced by Gromov \cite{Gromov} in the 1980s following earlier work of Trauber and Johnson), which admits powerful application in geometry, dynamics and rigidity theory (see e.g.\ \cite{Campagnolo-etal, Frigerio, Monod-Book}), but is notoriously hard to compute. The central conjecture concerning bounded cohomology of connected Lie groups is the
\begin{conj*}\label{conjecture} If $G$ is a connected semisimple Lie group with finite center, then the comparison map $c^\bcdot:\Hb^\bcdot(G) \to \HH^\bcdot(G)$ is an isomorphism in every degree. 
\end{conj*}
Here and in the sequel, given a topological group $G$, we denote by $\HH^{\bcdot}(G)$ the continuous group cohomology of $G$ with coefficients in the trivial $G$-module $\bbR$, and by $\Hb^\bcdot(G)$ the corresponding continuous bounded cohomology. In the present form, the isomorphism conjecture suggested by Monod in his ICM address in 2006 \cite{Monod-Survey}; surjectivity was already conjectured by Dupont \cite{Dupont}. For $G$ as in the isomorphism conjecture, the continuous cohomology $\HH^\bcdot(G)$ is isomorphic via the van Est isomorphism \cite[Cor. XI.5.6]{BW} to the cohomology of the compact dual of the non-compact symmetric space associated with $G$, which is known explicitly in all cases (see e.g. \cite{GHV,Toda-Mimura, Stas}). Combining this with the fact that continuous bounded cohomology is invariant under passing to the quotient by the amenable radical, the isomorphism conjecture makes an explicit prediction for the continuous bounded cohomology of all connected locally compact groups. 

In degrees $\leq 2$, the isomorphism conjecture was established by Burger and Monod in \cite{Burger-Monod1}. While surjectivity of the comparison map has been established for a large class of semisimple Lie groups, which includes for instance all groups of Hermitian type \cite{HartOtt2, BucherThesis}, proving injectivity in degrees $\geq 3$ is notoriously difficult, and there has been very limited progress over the last 15 years. 

For $\SL_2(\bbR)$ (respectively $\SL_2(\bbC)$), injectivity of the comparison map in degrees $3$ and $4$ (respectively degree $3$) has been established in \cite{Burger-Monod3} and \cite{HartOtt} (respectively in \cite{Burger-Monod3}, based on \cite{Bloch}). The methods extend to the reductive case (showing that  $ \Hb^3(\GL_2(\bbR)) \cong \{0\}$ and  $\Hb^3(\GL_2(\bbC)) \cong \bbR$), but not in any obvious way to the higher rank case. Nevertheless, injectivity in degree $3$ is also known for the groups $\SL_r(\bbR)$ and $\SL_r(\bbC)$ \cite{Monod-Stab, BBI}: the proof is by reduction to the rank one case using \emph{bounded-cohomological stability}. More precisely, it follows from stability results of \cite{Monod-Stab} that for all $r \geq 3$, we have injections
\[
\Hb^3(\SL_r(\bbR)) \hookrightarrow \Hb^3(\GL_2(\bbR)) = \{0\}  \qand \Hb^3(\SL_r(\bbC)) \cong \Hb^3(\GL_2(\bbC) \cong \bbR.
\]
Remarkably, prior to the current work, the aformentioned results were still the \emph{only}\footnote{It was pointed out by Monod that the results from \cite{Pieters} are not correct.}  injectivity results for the comparison map in degrees $\geq 3$ in the setting of the isomorphism conjecture. This hints at a close relation between bounded-cohomological stability and the isomorphism conjecture. In fact, this relation goes both ways:
\begin{enumerate}
\item The isomorphism conjecture \emph{predicts} that bounded-cohomological stability holds for all families of classical Lie groups for which continuous cohomology is stable. Besides general and special linear groups this concerns (general and special) orthogonal, unitary and symplectic groups (of increasing rank).
\item If bounded-cohomological stability (with a sufficiently good range) can be established for a family of Lie groups, then this can be \emph{used for computations} in continuous bounded cohomology and thereby help to establish various cases of the isomorphism conjecture, in particular in low degrees.
\end{enumerate}
In this context, the results of the present article and its two companion articles \cite{DMH0, DM} can be summarized as follows:
\begin{itemize}[leftmargin=25pt]
\item Bounded-cohomological stability can be established by a measurable version of Quillen's homological stability criterion; this is the main result of \cite{DMH0}.
\item Our version of Quillen's criterion can be applied to establish bounded-cohomological stability for many of the families, for which stability is predicted by the isomorphism conjecture; this is the goal of the present article.
\item Our stability results are strong enough to establish the isomorphism conjecture in degree $3$ for all classical complex connected semisimple Lie groups with finite center by reduction to the $\SL_2(\C)$ case; this is the main result of \cite{DM}.
\end{itemize}
\subsection{Main results}
We now explain the results and techniques of the present article in more details. Consider a sequence
\begin{equation} \label{eq:first_sequence}
G_0 < G_1 < G_2  < \dots
\end{equation} 
of Lie groups. We say that \eqref{eq:first_sequence} is \emph{bc-stable}\footnote{short for \emph{bounded-cohomologically stable}} if there exists a function $r: \bbN \to \bbN$ such that the respective inclusions induce isomorphisms
\[
\Hb^q(G_{r(q)}) \ \cong \ \Hb^q(G_{r(q)+1}) \ \cong \ \Hb^q(G_{r(q)+2}) \ \cong \ \cdots   
\]
in group homology for all $q \in \bbN$. Any such function $r$ is then called a \emph{stability range} for the family. As mentioned above, Monod has established that for $\kay \in \{\bbR, \bbC\}$, the family 
\[
\GL_1(\kay) < \GL_2(\kay) < \GL_3(\kay) < \cdots < \GL_{r}(\kay) < \cdots 
\]
is bc-stable with a linear stability range. As discussed in \cite{DMH0}, the same methods applies also to the quaternionic case $\kay= \bbH$, but not to any other classical families (due to the fact that only groups of type $A_n$ act highly transitively on products of generalized flag varieties). In the present article, we will develop new techniques, which will allow us to establish the following bounded-cohomological stability results, as predicted by the isomorphism conjecture.
\begin{thmx}\label{ThmIntro} 
The following sequences are bc-stable:
\begin{equation*}
\begin{array}{rll}
	\{1\} \!\!\!&< \Sp_2(\kay) < \Sp_4(\kay) < \cdots < \Sp_{2r}(\kay) < \cdots \quad &\mbox{for } \kay \in \{\bbR,\bbC\},\\[2pt]
	\OO_d(\bbC) \!\!\!& < \OO_{2+d}(\bbC) < \OO_{4+d}(\bbC) < \cdots < \OO_{2r+d}(\bbC) < \cdots &\mbox{for any } d \in \{0,1\},\\[2pt]
	\OO(d) \!\!\!& < \OO(d+1,1) < \OO(d+2,2) < \cdots < \OO(d+r,r)< \cdots \quad \, & \mbox{for any } d \in \bbN, \\[2pt]
	\UU(d) \!\!\!& < \UU(d+1,1) < \UU(d+2,2) < \cdots < \UU(d+r,r)< \cdots \quad &\mbox{for any } d \in \bbN.
\end{array}
\end{equation*}
\end{thmx}
In fact, the case of the symplectic groups was already established in an earlier preprint by the authors \cite{DM+Hartnick}. We present here an easier and more conceptual proof, which unlike the earlier proof generalizes to all of the families above.

From \autoref{ThmIntro}, we derive two corollaries. The first concerns the corresponding subgroups of matrices of determinant one:
\begin{corx} \label{ThmIntroSG}
The following sequences are bc-stable:
\begin{equation*}
\begin{array}{rll}
	\{1\} \!\!\! &< \SO_3(\bbC) < \SO_5(\bbC) < \cdots < \SO_{2r+1}(\bbC) < \cdots \\[2pt]
	\SO(d) \!\!\! &< \SO(d+1,1) < \SO(d+2,2) < \cdots < \SO(d+r,r)< \cdots  & \mbox{for \emph{odd} } d \in \bbN \\[2pt]
	\SO(d) \!\!\! &< \SO_0(d+1,1) < \SO_0(d+2,2) < \cdots < \SO_0(d+r,r)< \cdots  & \mbox{for \emph{odd} } d \in \bbN \\[2pt]
	\SU(d) \!\!\! &< \SU(d+1,1) < \SU(d+2,2) < \cdots < \SU(d+r,r)< \cdots  &\mbox{for any } d \in \bbN
\end{array}
\end{equation*}
\end{corx}
For $\kay \in \{\bbR, \bbC\}$, the corresponding statement for the family $(\SL_{r}(\kay))_r$ 
was established by Monod in \cite{Monod-Stab}; a similar arguments also applies in the quaternionic case $\kay = \bbH$ (see \cite{DMH0}).

Our second corollary concerns bounded-cohomological stability for families of lattices; for technical reasons this applies only in the case of connected \emph{simple} Lie groups.
\begin{rem} \label{simple_ones}
Among the ones listed in \autoref{ThmIntro} and \autoref{ThmIntroSG}, the following sequences consist of connected simple Lie groups (except for $r=0$):  
\begin{equation*} 
	\Sp_{2r}(\bbC), \ \ \Sp_{2r}(\bbR), \ \ \SO_{2r+1}(\bbC), \ \ \SO_0(d+r,r) \mbox{ for odd } d, \ \ \SU(d+r,r). \vspace{1pt}
\end{equation*}
\end{rem}
Combining \autoref{ThmIntro} and \autoref{ThmIntroSG} with \cite[Corollary 1.4]{Monod-sot} one obtains immediately the following consequence:
\begin{corx} \label{ThmLatticeIntro}
Let $(G_r)$ be one of the sequences in \autoref{simple_ones} and let $\Gamma_r < G_r$ be a lattice such that the inclusions $G_r \hookrightarrow G_{r+1}$ restrict to inclusions
\[
	\Gamma_0 < \Gamma_1 < \Gamma_2 < \cdots < \Gamma_r < \cdots
\]
Then $(\Gamma_r)_{r \in \bbN}$ is bc-stable.\qed
\end{corx}
 For example, the families $(\Sp_{2r}(\bbZ))_r$ and $(\Sp_{2r}(\bbZ[i]))_r$ are bc-stable. 

\subsection{Stiefel complexes associated with formed spaces}
We will give a unified proof for bc-stability of all of the families appearing in \autoref{ThmIntro} based on the fact that all of the
groups appearing in \autoref{ThmIntro} can be realized as automorphism groups of non-degenerate formed spaces. Here, by a (non-degenerate) \emph{formed space} $(V, \omega)$ over a field $\kay$ of characteristic $\neq 2$ we mean a finite-dimensional $\kay$-vector space $V$ together with a (non-degenerate) sesquilinear form $\omega$; the \emph{rank} of such a formed space is defined as the common dimension of all maximal totally $\omega$-isotropic subspaces of $V$. Our proof of  \autoref{ThmIntro} will be based on the following construction:
\begin{con}
We fix a non-degenerate formed space $(V, \omega)$ of rank $r$ with automorphism group $G = \Aut(V, \omega)$ over a field $\kay \in \{\bbR, \bbC\}$. We then denote by $\mathcal P = \mathcal P_{V, \omega}$ the space of $\omega$-isotropic lines in $V$, seen as a subset of the projective space $\bbP(V)$. We then define
\[
\bar X_{l}(V, \omega) := \{(p_0, \dots, p_l) \in \mathcal P_{V, \omega}  \mid \omega(p_i,p_j) = 0 \mbox{ for all } i , j \in \{0, \dots, l\}\} \quad (l \leq r-1),
\]
and denote by $X_l(V, \omega)$ the unique Zariski open $G$-orbit in $\bar X_l(V, \omega)$. 

Since $G$ acts transitively on $X_l(V, \omega)$ for every $l \leq r-1$, each of these spaces is a standard Borel space, which carries a unique $G$-invariant measure class. In the language of \cite{DMH0}, this means that each $X_l(V, \omega)$ is a \emph{Lebesgue $G$-space}. Moreover, for every $l \geq 1$, we have natural face maps $\delta_i: X_{l} \to X_{l-1}$ given by forgetting the $i$-th point, and the spaces $X_l$ together with these face maps form a semi-simplicial object $X(V, \omega)$ in the category of Lebesgue-$G$-spaces. In analogy with the orthogonal case \cite{Vogtmann}, we refer to $X(V, \omega)$ as the \emph{Stiefel complex} associated with $(V, \omega)$.
\end{con}
Every Stiefel complex $X = X(V, \omega)$ gives rise to a cochain complex of Banach spaces 
\begin{equation}\label{LInftyComplexStiefel}
0 \xrightarrow{\dd^{-2}} \bbR \xrightarrow{\dd^{-1}} L^\infty(X_0) \xrightarrow{\dd^{0}} L^\infty(X_1) \xrightarrow{\dd^{1}} L^\infty(X_2) \xrightarrow{\dd^{2}} L^\infty(X_3) \xrightarrow{} \cdots \xrightarrow{}L^\infty(X_{r-1}),
\end{equation}
called the \emph{augmented $\Linfty$-complex} of $X(V, \omega)$. Here, for every $l \leq r-1$, the Banach space $L^\infty(X_l)$ is defined using the invariant measure class on $X_l$, the maps $\dd^l: L^\infty(X_{l}) \to L^\infty(X_{l+1})$ are defined as the alternating sums of the dual face maps $\delta^i: L^\infty(X_{l}) \to L^\infty(X_{l+1})$, and $\dd^{-1}: \bbR \to  L^\infty(X_0)$ is the inclusion of constants. 
\begin{defn} The Stiefel complex $X(V, \omega)$ is \emph{boundedly $\gamma$-acyclic} for some $\gamma \leq r-1$, if the cochain complex \eqref{LInftyComplexStiefel} satisfies $\ker \dd_l = \im \dd_{l-1}$ for every $l \leq \gamma$. 
\end{defn}
Our main new tool in the current article will be the following acyclicity result:
\begin{thmx}\label{AcyclicityBound} If $(V, \omega)$ be a non-degenerate formed space of rank $r \geq 5$. Then the Stiefel complex $X(V, \omega)$ is $\gamma(r)$-acyclic, where
\begin{equation}\label{gammaIntro}
\gamma(r) \geq \sup \left\{l \mid  2^l + \left\lceil(l+1)/2\right\rceil \leq \left\lfloor \frac{r-1}{2}\right\rfloor\right\}.
\end{equation}
\end{thmx}
Note that the lower bound on $\gamma$ given by \eqref{gammaIntro} grows roughly like $\log_2(r)$ as $r \to \infty$. This is enough to ensure that $\gamma(r) \to \infty$ as $r \to \infty$, which is sufficient for our purposes. The acyclicity bound \eqref{gammaIntro} will be established by means of a random chain homotopy which is constructed by choosing points at random. This construction completely disregards the fine structure of the Stiefel complex, and hence the bound is certainly not sharp. We believe that a more careful analysis of the Stiefel complex could yield a linear acyclicity bound rather than the above logarithmic bound, but we will not pursue this here.

\subsection{From Stiefel acyclicity to stability}
Let $(G_r)_{r \geq 0}$ be one of the families from \autoref{ThmIntro}. For every $r \in \bbN$, there exists a non-degenerate formed space $(V_r, \omega_r)$ of rank $r$, such that $G_r = \Aut(V_r, \omega_r)$. Then $G_r$ acts $(r-1)$-transitively on the Stiefel complex $X(r) = X(V_r, \omega_r)$, and the actions are compatible in a suitable sense (see \autoref{thm:compatibility} below). In this situation the measurable Quillen criterion from \cite{DMH0} applies and states the following:
\begin{thmx}[{\cite[Theorem 4.6]{DMH0}}]\label{Quillen} Assume that for every $r \geq 0$, the Stiefel complex $X(r)$ is boundedly $\gamma(r)$-acyclic. If we define
\[
\widetilde{\gamma}(q, r) := \left\{\begin{array}{rl}\min \left\{\gamma\left(r+1-2(q-j)\right) - j \mid j \in \{1, \dots, q\}\right\}, & \text{if }r + 1 - 2(q-1) \geq 0,\\
-\infty, & \text{else},
 \end{array}\right.
\]
then the inclusion $\iota_r: G_r \to G_{r+1}$ induces an isomorphism (resp. an injection) 
\[\Hb^q(\iota_r): \Hb^q(G_{r+1}) \, \xrightarrow{\cong \ } \, \Hb^q(G_r) \quad (\text{resp. }\Hb^{q+1}(\iota_r): \Hb^{q+1}(G_{r+1}) \hookrightarrow \Hb^{q+1}(G_r)).
\] 
provided the pair $(q,r) \in \bbN^2$ satisfies the condition
\begin{equation}
\min\{\widetilde{\gamma}(q,r), r-2q\} \geq 0.
\end{equation}
In particular, if $\gamma$ is proper, then $\widetilde{\gamma}$ is proper and hence $(G_r)_{r \geq 0}$ is bc stable.
\end{thmx}
Now \autoref{ThmIntro} is immediate from \autoref{AcyclicityBound} and \autoref{Quillen}, and \autoref{ThmIntroSG} follows easily from \autoref{ThmIntro}.
\begin{rem} Because of our lazy acyclicity bounds in \autoref{AcyclicityBound}, the stability range obtained from  \autoref{AcyclicityBound} and \autoref{Quillen} is of the form $r(q) = \clO(2^q)$ as $q\to \infty$, whereas the isomorphism conjecture predicts a \emph{linear} stability range. One would expect that this causes a problem as far as computations of low degree continuous bounded cohomology using stability are concerned, but actually this is not the case: It will be proved in \cite{DM} that as far as concrete computations are concerned, \emph{any} stability range is as good as any other, due to a bootstrapping argument in continuous bounded cohomology. In fact, we are not aware of any applications of bc-stability which would benefit from a better range. For this reason, we  preferred to give a uniform proof of acyclicity which does not involve the fine structure of formed spaces rather than analyzing each case in detail to optimize the bounds.
\end{rem}
\begin{rem}
We have omitted from this article the treatment of the remaining quaternionic classical families $\Sp(d+r,r)$ and $\SO^\ast(2r)$. We expect our arguments to adapt readily to those groups after adapting the results from Subsection \ref{subsec:isotropic_pts}. 
\end{rem}
\subsection{Notational conventions} We convene that $0 \in \bbN$. For $k \in \bbN$, we denote by $[k]$ the set $\{0,\ldots,k\} \subset \bbN$, and set $[\infty] := \bbN$. We write $\Sym_k$ for the group of bijections $[k] \to [k]$. The curly letter $\kay$ will be reserved to denote either $\bbC$ or $\bbR$. We use the shorthands ``m.c.p.'' for ``measure-class preserving'', and ``p.m.p'' for ``probability-measure preserving'' maps. \vspace{5pt}

{\bf Acknowledgements.} We thank Marc Burger for suggesting the question that led to the main result of this paper and his input throughout the development of the project, as well as Uri Bader, Alessandra Iozzi, Stefan K\"uhnlein, Nir Lazarovich, Nicol\'as Matte Bon, Alexis Michelat, Nicolas Monod, Maria Beatrice Pozzetti, and Fabian Ruoff for insightful conversations. 

This work contains results of the first author's doctoral thesis \cite{DM-Thesis}, and was completed during his postdoctoral fellowship at the Weizmann Institute of Science, Rehovot, Israel. He was supported by the Swiss National Science Foundation, Grants No. 169106 and 188010. The second author was supported by the Deutsche Forschungsgemeinschaft, Grant No. HA 8094/1-1 within the Schwerpunktprogramm SPP 2026 (Geometry at Infinity).


\section{Formed spaces and their Stiefel complexes} \label{formed}
\subsection{Formed spaces} 
We recall some background concerning formed spaces, following \cite{Sprehn-Wahl1,Taylor}. We recall that $\kay$ denotes either $\bbC$ or $\bbR$.
In the sequel $\sigma \in \Aut(\kay)$ will denote either the identity or complex conjugation $\bar\cdot: \bbC \to \bbC$. 

Let $V$ be a finite-dimensional $\kay$-vector space, and let $\omega$ be a reflexive $\sigma$-sesquilinear form on $V$, that is, a function $\omega\!:\! V \times V \to \kay$ that is additive in each argument and satisfies \vspace{-2pt}
\[
	\omega(\alpha v, \beta w) = \sigma(\alpha) \,\beta \cdot \omega(v,w) \qqand
	\omega(v,w) = 0 \, \Longrightarrow \, \omega(w,v) = 0.
\]
for all $v, w \in V$ and $\alpha, \beta \in \kay$. Then the pair $(V, \omega)$ is called a \emph{formed space} (cf. \cite{Sprehn-Wahl1}). 

We say that $\omega$ is an \emph{isotropic} form if there exists a vector $v \in V \smallsetminus\{0\}$ such that $\omega(v,v) = 0$. For any $A \subset V$, let 
\begin{equation} \label{def_perp}
	A^\perp := \{v \in V \mid \omega(\{v\} \times A) = 0\}
\end{equation}
denote the $\omega$-\emph{perpendicular} of $A$. We say that a subspace $W < (V,\omega)$ is \emph{totally isotropic} if the inclusion $W \subset W^\perp$ holds, and \emph{non-degenerate} if $W \cap W^\perp = \{0\}$. If $(V,\omega)$ is non-degenerate, then $\dim(W) + \dim(W^\perp) = \dim(V)$ and hence $(W^\perp)^\perp = W$ for any $W < V$. We write $\mathrm{Rad}(\omega) := V^\perp$ for the \emph{radical} of $\omega$. If $V_0$ is any complementary subspace of the radical, then we have an $\omega$-orthogonal decomposition of the form $V = V_0 \oplus \mathrm{Rad}(\omega)$, with $\omega|_{V_0 \times V_0}$ non-degenerate.

An injective linear map $\varphi: (V_1,\omega_1) \to (V_2,\omega_2)$ between formed spaces is called a \emph{linear isometry} if $\varphi^\ast \omega_2 = \omega_1$. We write $G:=\Aut(V,\omega) < \GL(V)$ for the automorphism group of $(V,\omega)$, which consists of all bijective linear isometries $V \to V$. It is a closed subgroup of $\GL(V)$, hence a real Lie group with the subspace topology. Provided that $(V, \omega)$ is non-degenerate, the following version of Witt's lemma implies that all maximal (totally) isotropic subspaces of $(V,\omega)$ have the same dimension, called the \emph{rank} of $(V,\omega)$.

\begin{thm}[e.g. {\cite[Theorem 7.4]{Taylor}}] \label{thm:witt}
Let $(V,\omega)$ be a non-degenerate formed space, and $W < V$. Then any linear isometry $\varphi\!:\! (W,\omega\!\mid_{W \times W}) \!\to (V,\omega)$ extends to an element of $G$. \qed
\end{thm}

Let us assume that $(V,\omega)$ is non-degenerate. Then there exists $\varepsilon \in \kay$ with $\varepsilon \cdot \sigma(\varepsilon) = 1$ such that $\omega(w,v) = \varepsilon \, \sigma(\omega(v,w))$ for all $v,w \in V$; see \cite[Theorem 7.1]{Taylor}. We will record the parameters and call $(V,\omega)$ a \emph{$(\sigma,\varepsilon)$-formed space}. Up to rescaling $\omega$, we may assume that 
\[
(\sigma,\varepsilon)\in
	\{(\id_\kay,-1), \ (\id_\kay,+1), \ (\,\bar\cdot\,,+1) \},
\]
where the third case only occurs if $\kay = \C$. Correspondingly, we have three types of forms $\omega$: 
\begin{enumerate}[leftmargin=20pt]
\item $(\sigma,\varepsilon) = (\id_\kay,-1)$ means $\omega$ is \emph{alternating}, i.e. $\omega(v,v) =0$ for all $v \in V$; \vspace{1pt}
\item $(\sigma,\varepsilon) = (\id_\kay,+1)$ means $\omega$ is \emph{symmetric}, i.e. $\omega(v,w) = \omega(w,v)$ for all $v,w \in V$; or \vspace{1pt}
\item $(\sigma, \varepsilon) = (\, \bar\cdot\,,+1),\, \kay = \bbC,$ means $\omega$ is \emph{Hermitian}, i.e. $\omega(v,w) = \overline{\omega(w,v)}$ for all $v,w \in V$.
\end{enumerate}

\autoref{thm:adapted_basis} below states that the next two are the ``building blocks'' of any $(\sigma,\varepsilon)$-formed space, and we call them \emph{standard formed spaces}. 
\begin{enumerate}[leftmargin=20pt]
	\item \emph{$(\sigma,\varepsilon)$-hyperbolic space} $\clH_{\sigma,\varepsilon}$, which is the vector space $\kay^2$ endowed with the form \vspace{-4pt}
\[
	\omega^{\clH}_{\sigma,\varepsilon}(v,w) := \sigma(v_1) \, w_2 + \varepsilon \, \sigma(v_2) \, w_1, \quad \mbox{where } v=(v_1,v_2)^\top, \, w=(w_1,w_2)^\top.
\]
	\item \emph{$(\sigma,\varepsilon)$-Euclidean space} $\clE_{\sigma,\varepsilon}$, which is the vector space $\kay$, endowed with the form \vspace{-4pt}
	\[
		\omega^{\clE}_{\sigma,\varepsilon}(v,w) = (1+\varepsilon) \cdot \sigma(v)w,  \quad \mbox{where } v,w \in \kay.
	\] 

\end{enumerate}
Observe that $\clH_{\sigma,\varepsilon}$ has rank one, and that $\clE_{\sigma,\varepsilon}$ either is endowed with the identically zero form if $\varepsilon = -1$, or has rank zero if $\varepsilon = +1$. A \emph{hyperbolic pair} in a $(\sigma,\varepsilon)$-formed space $(V,\omega)$ is a pair of linearly independent vectors $(e,f)$ in $V$ that produces the coordinate-wise isomorphism $\Span\{e,f\} \cong \clH_{\sigma,\varepsilon}$. Analogously, a \emph{Euclidean vector} is a vector $h \neq 0$ that gives rise to a coordinate-wise isomorphism $\Span\{h\} \cong \clE_{\sigma,\varepsilon}$ of formed spaces. 

\begin{prop}[{e.g. \cite[\S 7]{Taylor}}] \label{thm:adapted_basis}
Let $(V,\omega)$ be a non-degenerate $(\sigma,\varepsilon)$-formed space of rank $r$ and dimension $n$, and $L=\Span\{e_1,\ldots,e_r\}$ be a maximal isotropic subspace. Then there are subspaces $L'=\Span\{f_1,\ldots,f_r\}$, $W=\Span\{h_1,\ldots,h_d\}$ of $V$, with $d := n-2r$, such that
\begin{enumerate}[label = \emph{(\roman*)},leftmargin=25pt]
	\item $L'$ is maximal isotropic and $W$ contains no isotropic vectors;
	\item all $(e_i,f_i)$ are hyperbolic pairs and all $h_j$ are Euclidean vectors;
	\item the spaces $\Span\{e_i,f_i\}$, $\Span\{h_j\}$ are pairwise $\omega$-perpendicular; and
	\item $V = L \oplus L' \oplus W = \bigoplus_i \Span\{e_i,f_i\} \oplus \bigoplus_j \Span\{h_j\}$. \qed
\end{enumerate}
\end{prop}

\begin{defn} \label{def:adapted_basis}
An ordered basis $\{e_r,\ldots,e_1,h_1,\ldots,h_d,f_1,\ldots,f_r\}$ as in \autoref{thm:adapted_basis} is called an \emph{adapted basis} of $V$. The matrix
\begin{equation*} \label{matrix_J}
J^{r,d}_\varepsilon:=\left(\begin{array}{ccc}
			 0 & 0 & Q_r\\[-2pt] 
			 0 & 1_d & 0 \\[-2pt]
			 \varepsilon Q_r & 0 & 0
			\end{array} \right) \in {\rm M}_{n}(\kay),
\end{equation*}
represents the form $\omega$ in such a basis, where $Q_r \in {\rm M}_{r}(\kay)$ is the matrix with 1's on its antidiagonal and zero elsewhere (the middle row and column are omitted if $d=0$). 
\end{defn}

\begin{rem} If $(V,\omega)$ is degenerate, we can choose a basis $\{o_1, \dots, o_l\}$ of $\mathrm{Rad}(\omega)$ and an adapted basis $\{e_r,\ldots,e_1,h_1,\ldots,h_d,f_1,\ldots,f_r\}$ of a non-degenerate complement $V_0$, and combine them into a basis of $V$. We refer to it as an \emph{adapted basis} of $(V, \omega)$ and to $r$ as the \emph{rank} of $(V, \omega)$ (note however that a maximal isotropic subspace of $V$ has dimension $r+l$).
\end{rem}
If $\varepsilon = +1$, we denote by $\clT_{\sigma, \varepsilon} := (\kay, 0)$ the 1-dimensional totally isotropic $(\sigma,\varepsilon)$-formed space over $\kay$. Then, we obtain:
\begin{cor}
Any $(\sigma,\varepsilon)$-formed space $(V,\omega)$ of rank $r$ and dimension $n$ is isomorphic to the orthogonal direct sum $\clH_{\sigma,\varepsilon}^{\oplus r} \oplus \clE_{\sigma,\varepsilon}^{\oplus d}\oplus \clT_{\sigma, \varepsilon}^{\oplus l}$, with $2r+d+l = n$, where the final summand appears only if $\varepsilon= +1$.
If $\varepsilon = +1$, $(V,\omega)$ is non-degenerate if and only if $l = 0$ and hence $d = n-2r$, and if $\varepsilon =-1$, it is non-degenerate if and only if $d=0$ and hence $n = 2r$.
\qed
\end{cor}

\begin{defn} We call $V_{\sigma,\varepsilon}^{r,d} := \clH_{\sigma,\varepsilon}^{\oplus r} \oplus \clE_{\sigma,\varepsilon}^{\oplus d}$ the (non-degenerate) \emph{standard} $(\sigma,\varepsilon)$-formed space of rank $r$ and dimension $2r+d$. The canonical basis of $V_{\sigma,\varepsilon}^{r,d}$ is an adapted basis. 
\end{defn}
\begin{rem} \label{restrictions_d}
In light of \autoref{thm:adapted_basis}, the field $\kay$ and the type $(\sigma,\varepsilon)$ of a non-degenerate formed space $(V,\omega)$ of rank $r$ and dimension $n$ impose additional restrictions on the natural number $d = n - 2r$ only in the following two cases:
\begin{enumerate}[leftmargin=25pt]
	\item If $(\sigma,\varepsilon) = (\id_\kay,-1)$, then $n = 2r$ is an even number.
	\item If $\kay = \bbC$ and $(\sigma,\varepsilon) = (\id_\bbC,+1)$, then $d \in \{0,1\}$. Indeed, two $\omega$-perpendicular Euclidean vectors $h_1, h_2 \in V$ would produce the isotropic vector $h_1 + ih_2$. 
\end{enumerate}
\end{rem}
\begin{rem}
All the aspects of the theory of formed spaces over fields above can be developed for any field $\kay$ of characteristic $\neq 2$ and any involution $\sigma \in \Aut(\kay)$, and extended to $\Char(\kay) = 2$ under consideration of an extra parameter; see \cite{Sprehn-Wahl1} for a succinct treatment in this generality. We assumed $\sigma$ to be continuous to avoid pathologies in the topology of $\Aut(V,\omega) < \GL(V)$. 
\end{rem}

\subsection{Notation} \label{subsec:grass}
Let $(V,\omega)$ be a formed space of positive rank $r$ and dimension $n$, and let $G=\Aut(V,\omega)$. We write $\bbP(V) := (V \smallsetminus \{0\})/\kay^\times$ for its projectivization, and $\bbP: V \smallsetminus \{0\} \to \bbP(V)$ for the quotient map. The image $\bbP(W)$ of a subspace $W<V$ is called a \emph{projective subspace} of $V$, and we call $\bbP(W)$ \emph{totally isotropic} if $W$ is totally isotropic with respect to $\omega$. We let
\begin{equation}
\clV := \{v \in V \mid \omega(v,v) = 0\} \qand \clP:= \bbP(\clV \smallsetminus \{0\}).
\end{equation}

If $k \in [n-1]$ and $W<V$ is a linear $\kay$-subspace of dimension $k+1$, then $\bbP(W)$ is a projective subspace of $\bbP(V)$ of dimension $k$, and we denote by $\Gr_k(V)$ the Grassmannian of all $k$-dimensional projective subspaces of $\bbP(V)$.\footnote{This is often denoted $\Gr_{k+1}(V)$ in the literature, but our difference in numbering will be advantageous later.} With our enumeration one has $\Gr_0(V) = \bbP(V)$ and $\Gr_{n-1}(V) = \{\bbP(V)\}$. We set $\Gr(V) := \Gr_0(V) \sqcup \cdots \sqcup \Gr_{n-2}(V)$. If $I \subset \bbP(V)$, let $\Span(I)$ be the point in $\Gr(V)$ corresponding to the $\kay$-subspace of $V$ spanned by $\bbP^{-1}(I) \subset V$. 

All the spaces $\Gr_k(V)$ for $k \in [n-1]$ are compact homogeneous $\GL(V)$-spaces, hence admit a continuous $G$-action by restriction. Given $k \in [r-1]$, we denote by $\clG_k \subset \Gr_k(V)$ the subset of all totally isotropic projective subspaces of dimension $k$, a compact homogeneous $G$-space.

\subsection{Isotropic points} \label{subsec:isotropic_pts}
The results about $\clV$ and $\clP$ presented here will be needed later.
\begin{lem}\label{linalg2} Let $(V, \omega)$ be a (possibly degenerate) $(\sigma, \varepsilon)$-formed space of rank $r$. Depending on $\kay$ and $(\sigma, \varepsilon)$, we make the following additional assumptions:
\begin{enumerate}[label=\emph{(\roman*)}]
\item If $\varepsilon = -1$, there is no assumption.
\item  If $\kay = \bbC$ and $(\sigma, \varepsilon) = (\overline{\cdot}, +1)$, we assume that $r \geq 1$.
\item If $\kay = \bbR$ and  $(\sigma, \varepsilon) = (\id_\bbR, +1)$, we assume that $r \geq 2$.
\item  If $\kay = \bbC$ and  $(\sigma, \varepsilon) = (\id_\bbC, +1)$, we assume that $r \geq 2$. \vspace{-4pt}
\end{enumerate}
Then $\Span(\clV) = V$ and, as a real projective algebraic set, $\clP$ is irreducible. Moreover, the assumptions are sharp, i.e.\ the conclusion may fail if $r = 0$ in (ii) or if $r = 1$ in (iii) and (iv).\vspace{-4pt}
\end{lem}
\begin{proof} If $(\sigma, \varepsilon) = (\id_\kay, -1)$, then $\omega$ is alternating, hence $\clV = V$, and the conclusion is evident.

Assume that $\kay = \bbR$ and $(\sigma, \varepsilon) = (\id_\bbR, +1)$. Then $\omega$ is a symmetric bilinear form, and the associated quadratic form $q$ can be diagonalized. Assuming positive rank, let $(e_1, \dots, e_n)$ be a basis of $V$ in which $q$ takes the form $q(x) = x_1^2 + \dots + x_r^2 - x_{r+1}^2 - \dots - x_{r+s}^2$ for $s\geq r \geq 1$ with $r+s \leq n$.
The $j$-th basis vector $e_j$ is then isotropic if $j \geq r+s+1$, contained in the hyperbolic plane $\Span\{e_j, e_{r+1}\}$ if $j \in \{1, \dots, r\}$, or in the hyperbolic plane $\Span\{e_1, e_{j}\}$ for $j \in \{r+1, \dots, r+s\}$. In particular, every basis vector $e_j$ is a linear combination of isotropic vectors, establishing $\Span(\clV) = V$. We consider now the algebraic set $\clP$, defined by the homogeneous polynomial $q$. It is irreducible if and only if the set given by intersection with the affine chart $\{x_1 = 1\}$ is irreducible. If $r=s=1$, this affine set is the zero set of $1-x_2^2= (1-x_2)(1+x_2)$, a reducible polynomial.
However, if $r \geq 2$, then it is the zero set of $p(x) = 1 + x_2^2 + \dots + x_r^2 - x_{r+1}^2 - \dots -  x_{r+s}^2$ for $s\geq r \geq 2$. We claim that the polynomial $p$ is irreducible. Indeed, assume that $p = p_1p_2$ for $p_1,\,p_2 \in \bbR[x_2, \dots, x_{r+s}]$. Then $1+x_2^2 = p_1(x_2, 0, \dots, 0) \cdot p_2(x_2, 0, \dots, 0)$, and since $1+x_2^2 \in \bbR[x_2]$ is irreducible, then one of the two factors, say $p_2$, must be independent of $x_2$. Writing $p_1(x) = Ax_2^2 + Bx_2 + C$ for some $A,B,C \in \bbR[x_3, \dots, x_{r+s}]$ and comparing coefficients of $x_2^2$ yield the equation $A\cdot p_2 = 1$, which shows that $p_2$ is a unit and, thus, that $p$ is irreducible.

Next, assume that $\kay = \C$ and $\varepsilon = +1$. Then $V$ admits a real form $V_\bbR$ such that $\omega|_{V_\bbR \times V_\bbR}$ is a bilinear form of the same rank $r>0$. Indeed, if $\sigma = \bar \cdot$, then we can choose $V_\bbR$ as the set of fixed points of the extension of $\sigma$ to $V$, and if $\sigma = \id_\bbC$, then we can choose $V_\bbR$ as the $\bbR$-span of any adapted basis. Then, by the previous case, we have $V_\bbR =  \Span_\bbR(V_\bbR \cap \clV) \subset \Span(\clV)$, and hence, $\Span(\clV)$ contains $V_\bbR \oplus i V_\bbR = V$.

We specify further to the case $\sigma = \, \bar \cdot$. In complex coordinates, $\clP$ is described by an equation of the form $|z_1|^2 + \dots + |z_{r}|^2 - |z_{r+1}|^2-\dots-|z_{r+s}|^2 = 0$ for some $s \geq r \geq 1$. In real coordinates, this amounts to $x_1^2 + y_1^2 + \dots + x_r^2 + y_r^2 - x_{r+1}^2-y_{r+1}^2 - \dots -x_{r+s}^2-y_{r+s}^2 = 0$. Hence, we can argue towards irreducibility as in the setting $\kay = \bbR$, $(\sigma,\varepsilon) = (\id_\bbR,+1)$.

Finally, we let $\sigma = \id_\bbC$. The set $\clP$ is described by an equation of the form $z_1^2 + \dots +z_{2r+d}^2 = 0$ for some $d \in \{0,1\}$. If we pass to real coordinates $x_j, y_j$ with $z_j = x_j + iy_j$, then the intersection of $\clP$ with the affine chart $\{z_1 = 1\}$ is described by the equation \vspace{-5pt}
\[
p(x_2, y_2, \dots, x_{2r+d}, y_{2r+d}) = (x_2y_2 + \dots + x_{2r+d}y_{2r+d})^2 - x_2^2-\dots- x_{2r+d}^2+y_2^2 + \dots + y_{2r+d}^2 -1 = 0.
\]

We need to investigate when $p$ is irreducible. If $r=1$ and $d = 0$ then $p(x_2, y_2) = (x_2^2+1)(y_2^2-1)$ is reducible. Thus, assume that $r \geq 2$, and suppose that $p$ factors as $p = qr$ with non-constant polynomials $q,r$. Write
\[
p = Ax_2^2 + Bx_2 + C \quad \text{with} \quad A,B,C \in R:=\bbR[y_2, \dots, x_{2r+d}, y_{2r+d}].
\]
We first note that $A(y_2, \dots, x_{2r+d}, y_{2r+d}) = y_2^2 -1 = (y_2+1) (y_2-1)$, but neither of these factors divides $C$, hence it is impossible that either $q$ or $r$ is constant in $x_2$.
The only other possibility is that $p = (Dx_2+E)(Fx_2+G)$ with $D,E,F,G  \in R$. This identity would then still hold if we set $y_3, x_4, y_4, \dots, x_{2r_d}, y_{2r_d}$ to $0$, so we would find $D', E', F', G' \in  \bbR[y_2, x_3]$ with
\[D'F' = y_2^2-1 = (y_2+1)(y_2-1), \quad D'G'+E'F' = 0, \quad E'G' =  -x_3^2-1.\]
Since $-x_3^2-1 \in \bbR[x_3]$ is irreducible, we may assume (up to reversing the order of the two factors and up to units) that $E' = x_3^2 + 1$ and $G' = -1$. Then $D' = E'F'$ is a multiple of $x_3^2+1$ which divides $y_2^2-1$, which is impossible. 
\end{proof}
\begin{lem}\label{linalg3} Let $(V, \omega)$ be a non-degenerate $(\sigma, \varepsilon)$-formed space of rank $r \geq 1$, and $W<V$ be a linear subspace of codimension $k$, which satisfy the following assumption:
\begin{enumerate}[label=\emph{(\roman*)}]
\item If $\varepsilon = -1$, there is no assumption.
\item If $\kay = \bbC$ and $(\sigma, \varepsilon) = (\overline{\cdot}, +1)$, then we assume that $k \leq r-1$.
\item If $\kay = \bbR$ and  $(\sigma, \varepsilon) = (\id_\bbR, +1)$, then we assume that $r \geq 2$ and $k \leq r-2$.
\item  If $\kay = \bbC$ and  $(\sigma, \varepsilon) = (\id_\bbC, +1)$, then we assume that $r \geq 3$ and $k \leq r - 2$.
\end{enumerate}
Then $\Span(W \cap \clV) = W$ and $\bbP\big((W \cap \clV) \smallsetminus\{0\}\big)$ is irreducible as a real algebraic set.
\end{lem}
\begin{proof}
The case $\varepsilon = -1$ is immediate, since $\clV = V$. In all other cases, we need to show in view of \autoref{linalg2} that the rank of $(W, \omega|_{W \times W})$ is at least $1$ or $2$, depending on $\kay$ and $(\sigma,\varepsilon)$. We may assume without loss of generality that $W$ is non-degenerate. Indeed: Choose an adapted basis $\{e_1,f_1,\ldots,e_r,f_r,h_1,\ldots,h_d\}$ of $V$ such that the radical $\mathrm{Rad}(W) = W \cap W^\perp$ is spanned by $e_1, \dots, e_a$. If $V_0$ is the subspace of $V$ spanned by $e_{a+1}, \dots, e_r, f_{a+1}, \dots, f_r, h_1, \dots, h_d$, then $V_0$ is non-degenerate and
\[
W \subset \mathrm{Rad}(W)^\perp = \mathrm{Rad}(W) \oplus V_0 \implies W = \mathrm{Rad}(W) \oplus (V_0\cap W).
\]
In particular, $W_0 := V_0 \cap W$ is a complement to the radical in $W$ and thus non-degenerate. If we denote by $k_0$ the codimension of $W_0$ in $V_0$, and by $ r_0= r-a$ the rank of $V_0$, we obtain
\[
r_0-k_0 = r_0-\dim_\kay V_0 + \dim_\kay W_0 = (r-a) - (\dim_\kay V-2a) + (\dim_\kay W - a) = r-k,
\]
and thus $k \leq r-j$ if and only if $k_0 \leq r_0-j$ for any $j$. Therefore, we may replace $V$ and $W$ by $V_0$ and $W_0$, respectively, and the difference $r-k$ will remain unchanged.

Our claim can then be established as follows, depending on the case at hand. First, we treat (ii) and (iii): Set $j:=1$ in case (ii) and $j:=2$ in case (iii). We have to show that $(W, \omega|_{W \times W})$ has rank at least $j$ provided $k \leq r-j$. Denote by $\sigma^\pm_V$ denote the number of positive, respectively negative eigenvalues of a matrix representing $\omega$ and define $\sigma^\pm_W$ accordingly. Since $V$ and $W$ are non-degenerate we have $\dim V =  \sigma^+_V+\sigma^-_V$ and  $\dim W = \sigma^+_W+\sigma^-_W$ and thus
\begin{align*}
\mathrm{rk}(W) &= \min\{\sigma^+_W, \sigma^-_W\} = \dim W - \max \{\sigma^+_W, \sigma^-_W\} \geq \dim W -  \max \{\sigma^+_V, \sigma^-_V\}
\\ & = \dim W - (\dim V - r) = r-k \geq r-(r-j) = j.
\end{align*}
For (iv), note that if $r \geq 3$ and $k \leq r-2$, then $k \leq 2r-5 \leq \dim_\C V - 5$, and hence $\dim_\C W \geq 5$. Since $W$ is non-degenerate, by Remark \ref{restrictions_d}, $W$ must have rank at least $2$.\end{proof}

\begin{lem} \label{linalg}
Let $k \in \bbN_0$ and let $p_0,\ldots,p_k,q_0,\ldots, q_k \in \clP$ be such that \vspace{-3pt}
\begin{equation} \label{linalg_eq}
	\forall \ i,j \in [k]: \quad \omega(p_i,p_j) = 0 \qand (\omega(p_i,q_j) = 0 \Longleftrightarrow i \neq j). \vspace{-3pt}
\end{equation}
Let $L < V$ be the linear subspace such that $\bbP(L) = \Span\{p_0,\ldots,p_k,q_0,\ldots,q_k\} \subset \bbP(V)$. Then $L$ is non-degenerate of dimension $2(k+1)$.  
\end{lem}
\begin{proof}
Our proof is by induction on $k$, where the base case $k=0$ is immediate. Our induction hypothesis is the statement for $k-1$. If $p_0,\ldots,p_k,q_0,\ldots,q_k$ are as in \eqref{linalg_eq}, then $L_0 < V$ defined by $\bbP(L_0):= \Span(p_0,\ldots,p_{k-1},q_0,\ldots,q_{k-1})$ is non-degenerate and of dimension $2k-1$. In particular, $V = L_0 \oplus L_0^\perp$. By assumption, $p_k \in \bbP(L_0^\perp)$. Now let $v_k \in p_k$ and $w_k \in q_k$ be non-zero vectors, and let $w_{k,0} \in L_0$ and $w_{k,1} \in L_0^\perp$ be such that $w_k = w_{k,0} + w_{k,1}$. Then 
$
	0 \neq \omega(v_k,w_k) = \omega(v_k, w_{k,1}), 
$
and hence, $\Span\{v_k,w_{k,1}\} \subset L_0^\perp$ is a $(\sigma,\varepsilon)$-hyperbolic space. Thus, $L=\Span(L_0 \cup  \{v_k,w_{k,1}\})$ is non-degenerate and of dimension $2(k+1)$. 
\end{proof}
\subsection{Stiefel complexes}
Let $(V,\omega)$ be a $(\sigma,\varepsilon)$-formed space over $\kay$, of rank $r > 0$ and dimension $n$, and let $G=\Aut(V,\omega)$. For any $l \in [r-1]$, we define
\begin{align*}
\BX_l &:= \{(p_0, \dots, p_l) \in \bbP(V)^{l+1}  \mid \omega(p_i,p_j) = 0 \mbox{ for all } i , j\}.
\end{align*}
Then the set $\BX_l$ is a real projective variety, thus a compact manifold.\footnote{If $\kay = \bbC$ and $\sigma = \id_\bbC$, it is acutally a complex projective variety, but we shall not rely on this structure.} 
Given $l \in [r-1]$ we say that points $p_0, \ldots, p_l \in \bbP(V)$ are \emph{in general position} if every subcollection $I \subset \{p_0,\ldots,p_l\}$ of points spans a subspace of $\bbP(V)$ of  dimension $\min\{|I|,r\}-1$. We then define a Zariski-open submanifold $X_l \subset \BX_l$ by
\begin{align*}
X_l &:= \{(p_0, \dots, p_l) \in \BX_l \mid p_0, \ldots, p_l \mbox{ are in general position}\}.
\end{align*}
It will be convenient to formally declare $\BX_{-1}$ and $X_{-1}$ to be singleton sets. Note that
$
	\BX_0 = X_0 = \clP
$
is the set of isotropic points in $\bbP(V)$. It equals $\bbP(V)$ if $\varepsilon = -1$ and is an embedded submanifold of $\bbP(V)$ of codimension $|\kay:\bbR|$ otherwise. 
\begin{defn} \label{def:Stiefel_variety}
The space $X_l$ is called the \emph{$l$-Stiefel manifold} of $(V,\omega)$, and $\BX_l$ is called the \emph{compactified} $l$-Stiefel manifold of $(V,\omega)$.  
\end{defn}
\begin{rem} \label{rem:bundle} For every $l \in \bbN$, the $l$-Stiefel manifold is a fiber bundle, $\Span: \, X_l \twoheadrightarrow \clG_s$, with $s := \min\{l,r-1\}$. The standard fiber is an open, dense subset of the product $\bbP(\kay^s)^{l+1}$. This justifies the terminology ``Stiefel manifold'' in analogy to the classical objects. 
\end{rem}
Both $\BX_l$ and $X_l$ admit $\Sym_{l}$-actions by permuting coordinates, and continuous $G$-actions by restricting the diagonal $G$-action on $\clP^{l+1}$. 
The following immediate consequence of \autoref{thm:witt} will be crucial for us:
\begin{lem} \label{thm:transitivity} For each $l \in [r-1]$, the automorphism group $G$ acts transitively on $X_l$. \qed
\end{lem}
Recall from \cite{DMH0} that a \emph{Lebesgue $G$-space} is a standard Borel space with a Borel action of $G$ and a $G$-invariant Borel measure class. As a consequence of \autoref{thm:transitivity}, each of the spaces $X_l$ with $l \in [r-1]$ is a Lebesgue $G$-space with respect to the unique $G$-invariant Borel measure class on $X_l$. Similarly, each of the spaces $\BX_l$ with $l \in [r-1]$ is a Lebesgue $G$-space with respect to the unique $G$-invariant measure class for which the subset $X_l$ is conull.

For $l \in [r-1]$ and $i \in [l+1]$, the maps 
\begin{equation} \label{eq:facemaps}
	\delta_i: \BX_{l+1} \to \BX_{l}, \quad \delta_i(p_0,\ldots,p_{l+1}) := (p_0,\ldots,\widehat{p_i},\ldots,p_{l+1})
\end{equation}
that delete the $i$-th coordinate are Borel, $G$-equivariant, and satisfy the face identities  
\begin{equation*} 
	\delta_{i} \, \circ \, \delta_{j} = \delta_{j-1} \, \circ \, \delta_{i} \qquad \mbox{whenever } i<j.
\end{equation*}
Their $G$-equivariance guarantees that they are measure class preserving (m.c.p.), and hence 
\[
\BX= ((\BX_l)_{l \in [r]}, (\delta_i: \BX_{l+1} \to \BX_{l})_{\underset{i \in [l+1]}{l \in [r-1]}}) \qand X= ((X_l)_{l \in [r]}, (\delta_i: X_{l+1} \to X_{l})_{\underset{i \in [l+1]}{l \in [r-1]}})
\]
are semi-simplicial objects in the category of Lebesgue $G$-spaces and $G$-equivariant m.c.p.\ Borel maps. Following \cite{DMH0} we refer to such objects as \emph{Lebesgue $G$-complexes}.
\begin{defn} The Lebesgue $G$-complex $X$ is called the \emph{Stiefel complex} of $(V,\omega)$.
\end{defn}
In the orthogonal case, this definition goes back to Vogtmann \cite{Vogtmann}. 


\section{Acyclicity of the Stiefel complex} \label{sec:contractibility}
\subsection{Statement of the main result}
The following notations and conventions will be used throughout this section: $(V,\omega)$ denotes a non-degenerate formed space of rank $r$ and dimension $n$. We set $d:= n-2r$ and $G := \Aut(V,\omega)$. Moreover, we define auxiliary functions
\begin{equation}
\begin{gathered}
\begin{array}{ll}
	\tilde\gamma_\ast(l) &:= 2^l + \lceil(l+1)/2\rceil, \ \ \ l \in \bbN, \\[2pt] 	
	\tilde\gamma\,(r)&:= \sup \gamma_\ast^{-1}\big((-\infty,r]\big), \ r \in \bbN. 
\end{array}
\end{gathered}
\end{equation}
Note that $\tilde\gamma_\ast: \bbN \to \bbN$ is a strictly increasing function, and that $\tilde\gamma_\ast(l) \leq r$ if and only if $l \leq \tilde\gamma(r)$. Moreover, $\tilde\gamma: \bbN \to \bbN \cup \{-\infty\}$ is increasing, yet not strictly. 
Now, set 
\[\gamma(r):= \tilde\gamma(\lfloor (r-1)/2 \rfloor) = \sup \left\{l \mid  2^l + \left\lceil(l+1)/2\right\rceil \leq \left\lfloor \frac{r-1}{2}\right\rfloor\right\}.
\]
Our goal is to establish \autoref{AcyclicityBound} from the introduction, which states that $X$ is boundedly $\gamma(r)$-acyclic. Since $X_l \subset \BX_l$ is conull, this is the same as proving that $\BX$ is boundedly $\gamma(r)$-acyclic, which is technically more convenient to prove. In fact, we will establish the following homotopy version:
\begin{thm}\label{thm:vanishing} Let $\gamma := \gamma(r)$. Then there exists a collection of bounded operators
\[
	h^l: \Linfty(\BX_{l+1}) \to \Linfty(\BX_l), \qquad l \in \{-1,0,\ldots,\gamma\}
\]
that satisfy the \emph{homotopy identity} $h^l \circ \dd^l + \dd^{l-1} \circ h^{l-1} = \id$ for all $l \in [\gamma]$.
\end{thm}
The remainder of this section is devoted to the proof of \autoref{thm:vanishing}. Since we will mostly working with $\BX$ instead of $X$, we will abuse language and also refer to $\BX$ as the Stiefel complex of $(V, \omega)$. From now on we work in the setting of \autoref{thm:vanishing} and denote $\gamma := \gamma(r)$.

\subsection{Heuristic for $h^l$} \label{heuristic}
Recall that each of the spaces $\BX_k$ comes equipped with a canonical $G$-invariant measure class and that $\BX_0 = \mathcal P$, the space of isotropic points. If we choose a $G$-quasi-invariant probability measure $\mu_0 \in \Prob(\BX_0)$ representing the invariant measure class on $\BX_0$, then $h^{-1}$ can be chosen as integration with respect to $\mu_0$. Based on this, we would like to construct inductively the maps $h^l$ for higher $l$. If $h^0$ is assumed to satisfy the homotopy identity for $l=0$, then \vspace{-2pt}
\begin{equation} \label{eq:heuristics1}
	h^0(\dd^0\!\varphi)(p_0) \overset{!}{=} \varphi(p_0) -\dd^{-1}h^{-1}\varphi(p_0) = \int_{\clP} (\varphi(p_0) - \varphi(t)) \, \dd\!\mu_0(t), \quad \varphi \in \Linfty(\BX_0). \vspace{-2pt}
\end{equation}
It is tempting to re-write the integrand as $\dd^0\!\varphi(t,p_0)$, which suggests formula
\[h^0f(p_0) =  \int_{\clP} f(t,p_0) \, \dd\!\mu_0(t) = \int_{\clP^2} f \, \dd(\mu_0 \otimes \mathrm{Dirac}_{p_0}).\] 
However, this re-writing is illegal since $f$ is only defined up to almost everywhere equivalence on $\mathcal P^2$, hence cannot be integrated against the singular measure $\mu_0 \otimes \mathrm{Dirac}_{p_0}$. The idea behind our solution is to rewrite
\[
\varphi(p_0) - \varphi(t) = (\varphi(s) - \varphi(t)) -  (\varphi(s) - \varphi(p_0)) = \dd^0\!\varphi(t,s) - \dd^0\!\varphi(p_0,s), 
\]
where $s \in \clP$ is an auxiliary random point. This random point should be chosen according to a suitable law with a Borel dependency on $(p_0,t) \in \clP^2$ such that
$
	\omega(p_0,s) = \omega(t,s) = 0
$
almost surely for $(\mu_0 \otimes \mu_0)$-almost every $(p_0,t) \in \clP^2$. Once we have established that such a law $\lambda_\perp^{p_0,t}$ exists, we may define
\begin{equation} \label{eq:heuristics2}
	h^0\!f(p_0) = \int_{\clP} \int_{\clP} \big(f(t,s)-f(p_0,s)\big) \dd\!\lambda_\perp^{p_0,t}(s) \, \dd\!\mu_0(t),
\end{equation}
and obtain \eqref{eq:heuristics1} as the special case $f := \dd^0\!\varphi$. 

Continuing the induction for higher values of $l$, one realizes the necessity of choosing more and more auxiliary, mutually dependent random points, and hence, of introducing some form of bookkeeping for all these dependencies. This will be the purpose of the next subsection.

\subsection{The configuration space $\scD_{l}$} \label{SecPBS}
Given a finite, totally ordered set $S$, we set $\Delta_0 S := S \sqcup 2^S$. On $\Delta_0 S$ we define an irreflexive, antisymmetric relation $<$ as follows:
\begin{itemize}[leftmargin=20pt]
	\item $<$ restricts on $S$ to the prescribed strict total order, and on $2^S$ to proper set inclusion;
	\item for $s \in S$ and $c \in 2^S$ we have $s < c$ if and only if $s \in c$.\end{itemize}
We call a subset $C \subset \Delta_0 S$ of cardinality $m+1$ an \emph{$m$-chain} if $<$ restricts to a total order on $C$; note that this implies in particular that $<$ is transitive on $C$ (whereas it is not transitive on $\Delta_0$!). For instance, if $S = \{0,1,2\}$ with the usual order, then $C := \{0 <  2 < \{0,2\} < \{0,1,2\}\}$ is a $3$-chain, but $\{0<1<2<\{0,2\}\}$ is not, since $1 \not < \{0,2\}$. The set $\Delta_0 S$ is canonically identified with the set of $0$-chains and we denote by $\Delta_m S$ the collection of all $m$-chains. Then the collection 
\[
	\Delta S := \bigsqcup_{m \in [|S|]} \Delta_m S
\] 
is an $|S|$-dimensional abstract simplicial complex with $m$-simplices being the elements of $\Delta_m S$. The orientation of simplices gives $\Delta S$ the structure of a semi-simplicial set in a natural way. Indeed, the $i$-th face operator
\begin{equation*} 
	\delta_i: \Delta_{m+1} S \to \Delta_m S \quad (i \in [m+1], \, \, m \in [|S| - 1]),
\end{equation*}
is the deletion of the $i$-th term in a chain. For instance, in the above example we have $\delta_1C = \{0<\{0,2\} < \{0,1,2\}\}$. In the sequel, we will use uppercase letters to denote chains, i.e.\ simplices of $\Delta S$, and reserve lowercase letters for vertices. Moreover, if $C$ is a simplex, we write $C_i$ for the $i$-th vertex in $C$, so that in the above example $C_2 = \{0,2\}$.

\begin{exmpl}\label{Deltal} If $S = [l-1]$, then $\Delta(l) := \Delta S$ is the  ``incomplete'' barycentric subdivision of an $l$-simplex in which one of its faces is kept undivided. For $l=1$ and $l=2$
 this is depicted in \autoref{fig:thepicture} and for $l=0$ we simply have $\Delta(0) = \{\emptyset\}$.
 \begin{center}
\begin{figure}[h!]
	\def\svgwidth{5.5cm}
	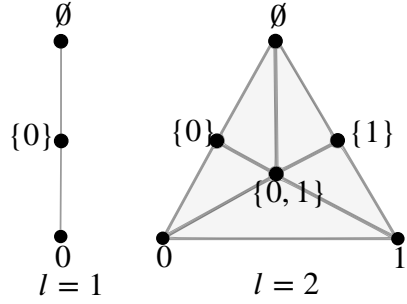 \vspace{-8pt}
	\caption{$\Delta(l)$ for $l= 1, 2$.}  \vspace{-2pt}
	\label{fig:thepicture} 
\end{figure} \vspace{-25pt}
\end{center} 
 \end{exmpl}

The definition of $\Delta$ is functorial in the following sense: Any order-preserving map $f: T \to S$ extends to a function $\Delta_0T \to \Delta_0S$ such that the induced relations are preserved. This enables the vertex-wise definition of an (injective) simplicial map $\Delta f: \Delta T \to \Delta S$, hence automatically a semi-simplicial morphism. 

\begin{defn} \label{def:config_space} For $l \in \bbN$ we define the \emph{configuration space} $\scD_{l} := \Hom(\Delta(l),\BX)$ as the set of morphisms of $\Delta(l)$ into the Stiefel complex $\BX$ (regarded as a semi-simplicial set), where $\Delta(l)$ is as in \autoref{Deltal}.
\end{defn}

Observe that $\scD_l \neq \emptyset$ for it contains the morphism that is constantly equal to $(p_0,\ldots,p_0) \in \BX_m$ on each $\Delta_m(l)$, $m\leq l,$ for any fixed $p_0 \in \clP$. We endow $\scD_l$ with the diagonal $G$-action, and topologize it as a closed (hence compact) subspace of the compact space $\BX^{\Delta(l)} := \prod_{m = 0}^l (\BX_m)^{\Delta_m(l)}$. If $t \in \scD_l$ and $C \in \Delta_m(l)$, we will write $t_C := t(C)$, and let \vspace{-2pt}
\begin{equation} \label{eq:eval}
	\ev_C: \scD_l \to \BX_m, \quad \ev_C(t):= t_C \vspace{-3pt}
\end{equation}
be the evaluation map, continuous by definition.

\begin{exmpl} \label{exmpl_delta}
Since $\Delta(0) = \{\emptyset\}$ we have $\scD_0 \cong \clP$. For $l = 1$ we read off from \autoref{fig:thepicture} that the collection of 0- and 1-simplices of $\Delta(1)$ are \vspace{-2pt}
\[
  	\Delta_0(1) = \{0, \emptyset, \{0\}\} \qand \Delta_1(1) = \{(0 < \{0\}), (\emptyset < \{0\})\}, \vspace{-2pt}
\]
respectively, hence $\scD_1$ will be the configuration space of all vertices in $\BX$ subjected to the constraints given by the edges in $\Delta(1)$. One verifies that $\scD_1$ is in bijection with the set\vspace{-2pt}
\begin{equation} \label{points_delta1}
	\{ (t_0, t_\emptyset, t_{\{0\}}) \in \clP^3 \mid \omega(t_{0},t_{\{0\}}) = \omega(t_{\emptyset},t_{\{0\}}) = 0\}. \vspace{-2pt}
\end{equation}
These are precisely the conditions imposed on the parameters $(p_0,t,s)$ in Formula \eqref{eq:heuristics2} from our heuristic discussion! 
\end{exmpl}

\begin{rem} \label{rem:incarnations}
There is an alternative way to think of the configuration space $\scD_l$ which will become important later. By definition, every element of $\scD_l$ is a simplicial map $t$ from $\Delta(l)$ to $\BX$. Restricting to vertices we thus obtain a map $\mathrm{res}(t) := t|_{\clP^{\Delta_0(l)}} :  \Delta_0(l) \to \clP = \BX_0$. Since $t$ is simplicial, the restriction $\mathrm{res}(t)$ is contained in the subset
\begin{equation}\label{scS}
\scS(l) := \big\{t\!:\! \Delta_0(l) \to \clP \mid \forall \,  (b < c) \in \Delta_1(l) : \, (t_b, t_c) \in \BX_1\big\} \subset \clP^{\Delta_0(l)}.
\end{equation}
We now observe that the map $\mathrm{res}: \scD_l \to \scS(l)$ is a $G$-equivariant homeomorphism, which allows us to identify $\scD_l$ with $\scS(l)$. Indeed, an element $(p_0,\ldots,p_k) \in \clP^{k+1}$ belongs to $\BX_k$ if and only if $(p_i,p_j) \in \BX_1$ for every $i, j \in [k]$, and this allows us to extend every map from the vertices of $\scD_l$ onto $\scS(l)$ to a unique simplicial map from $\scD_l$ to $\BX$.
\end{rem}

The following discussion concerns the shape of $\Delta S$ and will be useful in \autoref{interseccion}. 

\begin{rem} \label{rem:star}
Every codimension-one simplex of $\Delta S$ is the face of at most two top-dimensional simplices. More precisely, let $l := |S|$ and $C \in \Delta_l S$. There then exist $k \in \{-1,\ldots,l\}$ and elements 
$s_0,\ldots,s_l \in S$ and chains $C_{k+1}, \dots, C_l$ such that
 such that \vspace{-2pt}
\begin{equation}\label{TypicalC}
	C= (s_0 < \cdots < s_k < C_{k+1} < \cdots < C_l) \in \Delta_l S, \vspace{-2pt}
\end{equation}
with $C_{k+1} = \{s_0,\ldots,s_k\}$ and $C_{i+1} \smallsetminus C_i = \{s_i\}$ for $i \geq k+1$. Then $C$ is the only top simplex in $\Delta S$ with $\delta_l(C)$ as a face. Moreover, if $j \in [l-1]$ then there exists a unique simplex $C' \in \Delta_l S \setminus\{C\}$ and a unique $m=m(j) \in [l-1]$ such that
\[
\delta_j(C) = \delta_{m(j)}(C').
\]
Indeed, depending on $j$ and $k$, the pair $(C', m)$ is given as follows:
\begin{enumerate}[leftmargin=60pt]
	\item[$(j \!<\! k+1) \ $] $C' \!:= (s_0 <\!\cdots\! < \widehat{s_j} < \!\cdots\! < s_k < \{s_0, \ldots, \widehat{s_j}, \ldots, s_k \} < C_{k+1} < \!\cdots\! < C_l), \ \ m := k$;
	\item[$(j \!=\! k+1) \ $] $C' \!:= (s_0 < \!\cdots\! < s_k < s_{k+1} < C_{k+2} < \!\cdots\! < C_l), \ \ m := j = k+1;$
	\item[$(j \!>\! k+1) \ $] $C' \!:= (s_0 < \!\cdots\! < s_k < C_{k+1} < \!\cdots\! < C_{j-1} \!< C_{j-1} \cup \{s_j\} \!< C_{j+1} <\!\cdots\!< C_l), \ m := j$. 
\end{enumerate}
In each of the cases, the vertices $C_j$ and $C'_{m(j)}$, which are not in the shared face, are not joined by an edge in $\Delta S$. 
\end{rem}
\begin{defn} \label{def:star}
Let $l := |S|$ and $C \in \Delta_l S$ as in \eqref{TypicalC}. Given $C' \in \Delta_l S$ we write $C \sim C'$ if $C'$ is equal to $C$ or shares a codimension $1$ face with $C$. Then the \emph{star} $\St(C)$ of $C$ and the \emph{relative star} $\St(B)$ of $B := \delta_l(C) \in \Delta_{l-1} S$ are defined as \vspace{-3pt}
\[
\St(C) := \bigcup_{C' \sim C} \bigcup_{i=1}^l C'_i \quad \text{and} \quad\St(B) := \St(C) \smallsetminus \{S\}. 
\]
Since for every $B \in \Delta_{l-1}S$ there is a unique $C \in \Delta_lS$ such that $B = \delta_l(C)$ this is well-defined. Moreover, if $j \in [l-1]$, then by \autoref{rem:star} there is a unique $C' \in \Delta_l S$ and a unique $m \in [l-1]$ such that $\delta_j(C) = \delta_{m}(C')$, and we set
$
\mathrm{op}_j(C) := \mathrm{op}_j(B) := C'_m,
$
the \emph{vertex opposite $C_j$ with respect to $C$}. 
\end{defn}
\begin{rem}\label{StarBExplicit} By definition, we have \vspace{-3pt}
 \[
\St(C) = \{C_l\} \cup \bigcup_{i=0}^{l-1} \{C_i, \mathrm{op}_i(C)\} \qand \St(B) = \bigcup_{i=1}^{l-1} \{B_i, \mathrm{op}_i(B)\}.
\]
In particular, $|\St(C)| = 2l+1$ is always odd and $|\St(B)|=2l$ is always even.
\end{rem}
\begin{exmpl} Let $C:=(0\!<\!\{0\}\!<\!\{0,1\}\!<\!\{0,1,2\}) \in \Delta_3(2)$ so that $B =\delta_3(C) = (0\!<\!\{0\}\!<\!\{0,1\})$. Then
\[
\St(C) = \{0, \emptyset, \{0\}, 1, \{0,1\}, \{0,2\} \{0,1,2\}\} \qand \St(B) = \{0, \emptyset, \{0\}, 1, \{0,1\}, \{0,2\}\}.
\]
In this list, every vertex is followed by its opposite. \autoref{fig:theotherpicture} visualizes this example; the blue line goes over all the vertices contained in $\St(B)$.
\begin{center}
\begin{figure}[h!]
	\def\svgwidth{8cm}
\begingroup%
  \makeatletter%
  \providecommand\color[2][]{%
    \errmessage{(Inkscape) Color is used for the text in Inkscape, but the package 'color.sty' is not loaded}%
    \renewcommand\color[2][]{}%
  }%
  \providecommand\transparent[1]{%
    \errmessage{(Inkscape) Transparency is used (non-zero) for the text in Inkscape, but the package 'transparent.sty' is not loaded}%
    \renewcommand\transparent[1]{}%
  }%
  \providecommand\rotatebox[2]{#2}%
  \newcommand*\fsize{\dimexpr\f@size pt\relax}%
  \newcommand*\lineheight[1]{\fontsize{\fsize}{#1\fsize}\selectfont}%
  \ifx\svgwidth\undefined%
    \setlength{\unitlength}{985.90049072bp}%
    \ifx\svgscale\undefined%
      \relax%
    \else%
      \setlength{\unitlength}{\unitlength * \real{\svgscale}}%
    \fi%
  \else%
    \setlength{\unitlength}{\svgwidth}%
  \fi%
  \global\let\svgwidth\undefined%
  \global\let\svgscale\undefined%
  \makeatother%
  \begin{picture}(1,0.52018977)%
    \lineheight{1}%
    \setlength\tabcolsep{0pt}%
    \put(0,0){\includegraphics[width=\unitlength,page=1]{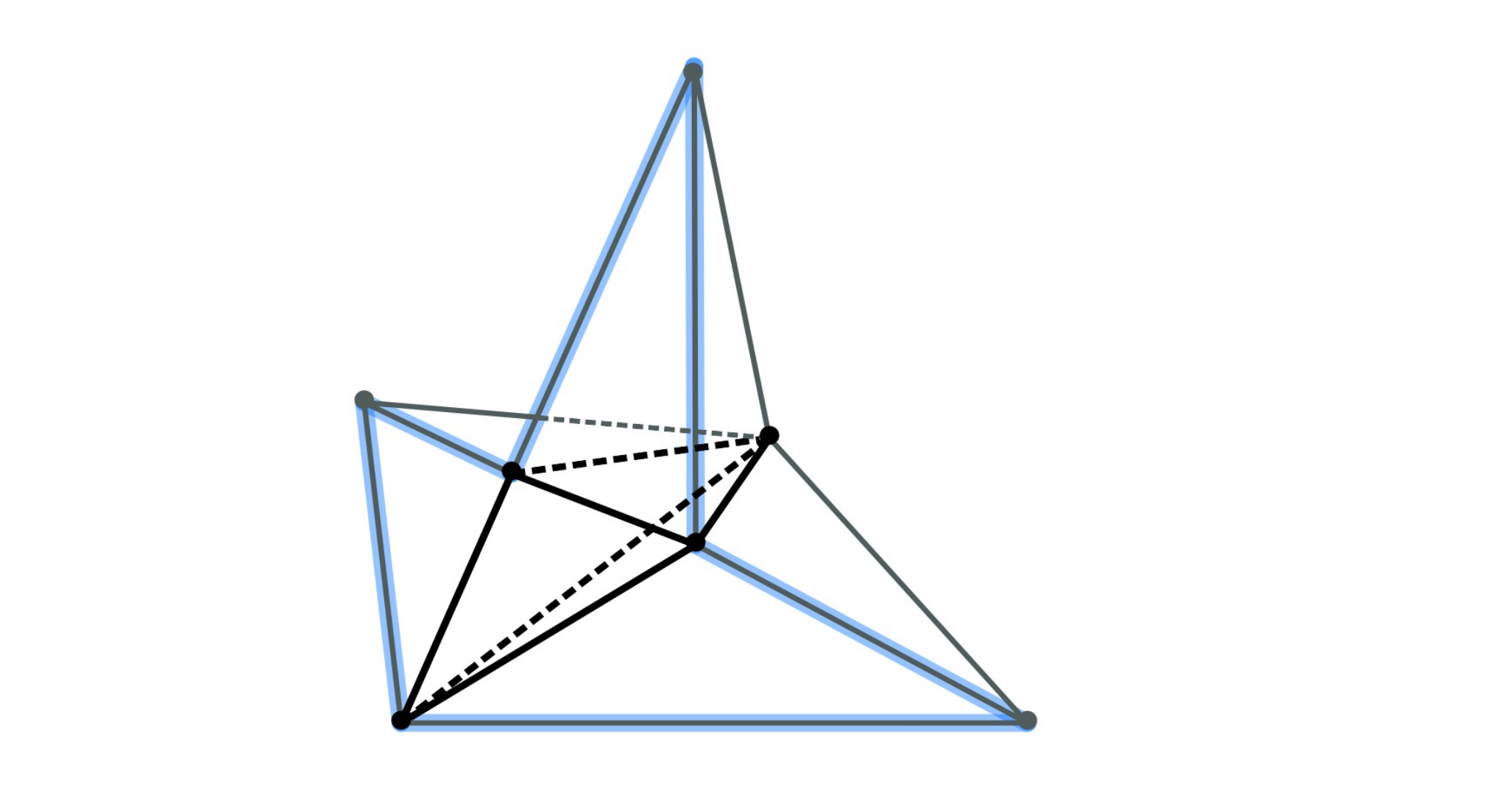}}%
    \put(0.48332076,0.45951225){\color[rgb]{0,0,0}\makebox(0,0)[lt]{\lineheight{1.25}\smash{\begin{tabular}[t]{l}$\mathrm{op}_0(B)=\emptyset$\end{tabular}}}}%
    \put(0.099915,0.03076836){\makebox(0,0)[lt]{\lineheight{1.25}\smash{\begin{tabular}[t]{l}$B_0 = 0$\end{tabular}}}}%
    \put(0.1331404,0.17239052){\makebox(0,0)[lt]{\lineheight{1.25}\smash{\begin{tabular}[t]{l}$B_1=\{0\}$\end{tabular}}}}%
    \put(0.44385924,0.10551715){\makebox(0,0)[lt]{\lineheight{1.25}\smash{\begin{tabular}[t]{l}$B_2 =\{0,1\}$\end{tabular}}}}%
    \put(0.696171,0.03532951){\makebox(0,0)[lt]{\lineheight{1.25}\smash{\begin{tabular}[t]{l}$\mathrm{op}_1(B)=1$\end{tabular}}}}%
    \put(0.51771927,0.23137942){\makebox(0,0)[lt]{\lineheight{1.25}\smash{\begin{tabular}[t]{l}$\{0,1,2\}$\end{tabular}}}}%
    \put(-0.00139664,0.27819106){\makebox(0,0)[lt]{\lineheight{1.25}\smash{\begin{tabular}[t]{l}$\mathrm{op}_2(B) = \{0,2\}$\end{tabular}}}}%
    \put(0.3538023,0.13573781){\makebox(0,0)[lt]{\lineheight{1.25}\smash{\begin{tabular}[t]{l}$B$\end{tabular}}}}%
  \end{picture}%
\endgroup%
 \vspace{-10pt}
	\caption{$\St(B)$ for $B = (0 < \{0\} < \{0,1\}) \in \Delta(2)$.} 
	\label{fig:theotherpicture}
\end{figure}
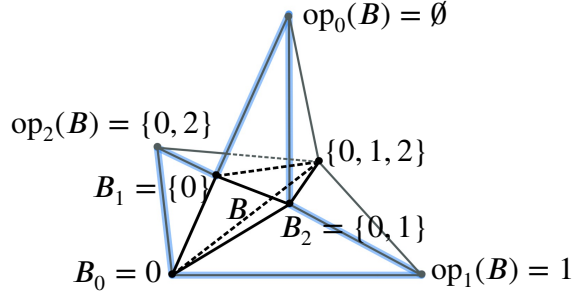 \vspace{-20pt}
\end{center} 
\end{exmpl}
For any $l \in \bbN$, we now define continuous $G$-equivariant maps
\begin{equation}\label{eq:partial}
\begin{gathered}
\begin{array}{ll}
	\partial_\emptyset:& \, \scD_{l+1} \to \BX_l, \quad \partial_\emptyset := \ev_{(0<\cdots<l)}, \qand\\[2pt] 
	\partial_i:& \, \scD_{l+1} \to \scD_{l}, \quad \partial_i t := t \, \circ \, \iota_i \quad (i \in [l]),
\end{array}
\end{gathered}
\end{equation}
where $\iota_i: \Delta(l) \rightarrow \Delta(l+1)$ is the injection induced by the order-preserving inclusion $[l] \smallsetminus \{i\} \hookrightarrow [l]$, after identifying $\Delta(l) \cong \Delta([l] \smallsetminus \{i\})$ via functoriality. They fit in the commutative square:
\begin{equation} \label{eq:interaction}
\begin{gathered}\xymatrix@C+3pc@R-6pt{
\scD_{l+1} \ar[r]^{\partial_i} \ar[d]_{\partial_\emptyset}& \scD_l \ar[d]^{\partial_\emptyset}\\
\BX_l \ar[r]^{\delta_i}& \BX_{l-1}
}\end{gathered}
\end{equation}
Our next goal is to equip these spaces with suitable probability measures which are compatible with face maps in a suitable sense. The following is our wishlist concerning such measures: 
\begin{thm}\label{PBSMeasures}
There exist $G$-quasi-invariant Borel probability measures $\mu_l \in \Prob(\BX_l)$ for all $l \in[r]$ and Borel probability measures $\sigma_l \in \Prob(\scD_l)$ for all $l \in [\gamma+1],$ such that:
\begin{enumerate}[label=\emph{(\roman*)},leftmargin=25pt]
	\item For every $l \in [r]$, the measure $\mu_l$ assigns full measure to $X_l$.
	\item For every $l \in [r]$, the face maps $\delta_i: (\BX_{l}, \mu_l) \to (\BX_{l-1}, \mu_{l-1})$ are p.m.p.
	\item For every $l \in [\gamma+1]$ and $C \in \Delta_l(l)$, the map $\ev_C: (\scD_l, \sigma_l) \to (\BX_l, \mu_l)$ is m.c.p.
	\item If $l \in [\gamma]$, then $\partial_\emptyset: (\scD_{l+1},\sigma_{l+1}) \to (\BX_l, \mu_l)$ is p.m.p.
	\item If $l \in [\gamma]$, then $\partial_i: (\scD_{l+1},\sigma_{l+1}) \to (\scD_l,\sigma_l)$ is p.m.p. for every $i \in [l]$.
\end{enumerate}
\end{thm}
We will demonstrate in Subsections \ref{subsec:perps} and \ref{SecMeasuresPBS} by means of an explicit construction that this wishlist can be fulfilled. However, it is important not to ask the evaluation maps in (iii) to be p.m.p., since we cannot ensure the existence of the measures $\sigma_l$ otherwise. For the moment, let us take the measures $\mu_l$ and $\sigma_l$ for granted and proceed to prove \autoref{thm:vanishing} along the lines suggested by our heuristic. Relying on \autoref{PBSMeasures} (iv), for any $l \in [\gamma]$, let \vspace{-1pt}
\begin{equation} \label{sigma_disint}
	\sigma^{(-)}: \BX_l \to \Prob(\scD_{l+1}) \vspace{-2pt}
\end{equation}
be the \emph{disintegration} of the measure $\sigma_{l+1}$ along the fibers of the p.m.p. map $\partial_{\emptyset}$, i.e. the unique Borel map (modulo $\mu_l$-null sets) such that
\begin{enumerate}[label=(\roman*),leftmargin=25pt]
	\item the fiber $\partial_\emptyset^{-1}(\mathbf{p})$ has full $\sigma^{\mathbf{p}}$-measure for $\mu_l$-a.e. $\mathbf{p} \in \BX_l$, and
	\item $\sigma_{l+1} = \int_{\BX_l} \sigma^{\mathbf{p}} \, \dd\!\mu_l(\mathbf{p})$, meaning that for every $f \in \scL^1(\scD_{l+1})$, \vspace{-1pt}
	\begin{equation} \label{defdisint}
		\int_{\scD_{l+1}} f \, \dd\!\sigma_{l+1} = \int_{\BX_l} \left( \int_{\scD_{l+1}} f \, \dd\!\sigma^{\mathbf{p}}\right)\dd\!\mu_l(\mathbf{p}).
	\end{equation}
\end{enumerate}
Recall that for any compact, separable, metrizable space $X$, the space $\Prob(X)$ has a Borel structure, regarded as a compact subspace of the dual Banach space $C(X)^\ast$, endowed with the weak-$\ast$ topology. For a reference on measure disintegration, see e.g. \cite[\S 5.3]{EW}.  

The fact that the maps in \autoref{PBSMeasures} (iv) and (v) are p.m.p. will be applied in the form of the next corollary in the proof of \autoref{thm:vanishing}, at the end of Subsection \ref{homotopy}.

\begin{cor} \label{partial_fiber}
For any $l \in [\gamma]$ and any $i \in [l]$, the equality $\partial_{i \, \ast} \sigma^{\mathbf{p}} = \sigma^{\delta_i(\mathbf{p})}$ of fiber measures on $\scD_l$ holds for almost all $\mathbf{p} \in \BX_l$. 
\end{cor}
\begin{proof}
All arrows in the commutative diagram \eqref{eq:interaction} are p.m.p. by Parts (ii), (iv) and (v) of \autoref{PBSMeasures}, hence the naturality statement \cite[Corollary 5.24]{EW} applies to the map $\partial_i$. 
\end{proof}

\subsection{The subdivision map $\beta^l$ and the proof of \autoref{thm:vanishing}} \label{homotopy}
For this section, we assume that probability measures $\mu_l$ and $\sigma_l$ as in \autoref{PBSMeasures} are fixed. Let $l \in \{-1,\ldots,\gamma\}$, and let $\Phi \subset \Delta_{l+1}(l+1)$ be the subset of all top simplices whose least element is the empty set. The action of the symmetric group $\Sym_{l}$ on $2^{[l]}$ induces a simply-transitive action on $\Phi$, hence for every $C \in \Phi$ there is a unique permutation $\alpha_C \in \Sym_{l}$ such that \[\alpha_C \cdot C = (\emptyset < [0] < [1] < \cdots < [l]).\]
\begin{defn}
Let $C \in \Delta_{l+1}(l+1)$; write $C$ as 
\begin{equation} \label{stdform}
	C= (s_0 < \cdots < s_k < C_{k+1} < \cdots < C_{l+1}), \vspace{-2pt}
\end{equation}
with $k \in \{-1,\ldots,l\}$ elements $s_0,\ldots,s_{l+1} \in [l]$ and $C_{k+1} = \{s_0,\ldots,s_k\}$. Moreover, assume that $C_{i+1} \smallsetminus C_i = \{s_i\}$ for all $i \geq k+1$. We then define $\phi(C) \in \Phi$ as
\[
	\phi(C) := (\emptyset < \{s_0\} < \cdots < \{s_0,\ldots,s_{k-1}\} < C_{k+1} < \cdots < C_{l+1}) \in \Phi,
\]
and set $\sgn(C):= (-1)^{k+1} \cdot \sgn(\alpha_{\phi(C)})$, where $\sgn(\alpha_{\phi(C)})$ is the signature of the permutation $\alpha_{\phi(C)} \in \Sym_{l}$. Finally, we define the \emph{subdivision map} 
\[
	\beta^{l+1}: \Linfty(\BX_{l+1}) \to \Linfty(\scD_{l+1}), \quad \beta^{l+1} := \sum_{C \in \Delta_{l+1}(l+1)} \sgn(C) \cdot \Linfty(\ev_C), \vspace{-4pt}
\]
where $\Linfty(\ev_C): \Linfty(\BX_{l+1}) \to \Linfty(\scD_{l+1})$ is the map induced by $\ev_C: \scD_{l+1} \to \BX_{l+1}$, which is well-defined because of \autoref{PBSMeasures} (iii). 
\end{defn}

\begin{exmpl} \label{exmpl_beta}
If $l = 0$ above, then $\beta^1\!f(t) = f(t_\emptyset,t_{\{0\}}) - f(t_0,t_{\{0\}})$ for $f \in \Linfty(\BX_1)$ and $t \in \scD_1$. Up to relabelling, this is the integrand of the expression \eqref{eq:heuristics2} from the heuristic discussion. Analogously, for $l = 1$, we have \vspace{-2pt}
\[
	\beta^1\!f(t) = f(t_0,t_1,t_{\{0,1\}}) - f(t_0,t_{\{0\}},t_{\{0,1\}}) + f(t_\emptyset,t_{\{0\}},t_{\{0,1\}}) - f(t_\emptyset,t_{\{1\}},t_{\{0,1\}}) + f(t_1,t_{\{1\}},t_{\{0,1\}}). \vspace{2pt}
\]
\end{exmpl}

Given $j \in [l-1]$, let us denote by $\partial^\emptyset\!:\! \Linfty(\BX_l) \to \Linfty(\scD_{l+1})$ and $\partial^j\!:\! \Linfty(\scD_l) \to \Linfty(\scD_{l+1})$ the maps at the level of $\Linfty$-spaces, induced by the maps $\partial_\emptyset$ and $\partial_j$ from \eqref{eq:partial}. These are well defined by \autoref{PBSMeasures} (iv) and (v). Now, for $l \in [\gamma]$, consider the diagram \vspace{-5pt}
\begin{equation*} 
\begin{gathered}\xymatrix@C+4pc{
\Linfty(\BX_l) \ar[r]^{\dd^l} \ar[rd]^{\partial^\emptyset} \ar[d]_{\beta^l}& \Linfty(\BX_{l+1}) \ar[d]^{\beta^{l+1}}\\
\Linfty(\scD_l) \ar[r]^{\sum_{j=0}^{l} (-1)^j \partial^j}& \Linfty(\scD_{l+1})
}\end{gathered}
\end{equation*}

\begin{lem}\label{thm:combinatorial}
For all $l \in [\gamma]$, the equation $\beta^{l+1} \circ \dd^l + \sum_{j=0}^l (-1)^j \partial^j \circ \beta^l = \partial^\emptyset$ holds. 
\end{lem}
\begin{proof}
Let $f \in \Linfty(\BX_l)$ and $t \in \scD_{l+1}$. Then the equality
\begin{equation}\label{computationhomot}
\big(\beta^{l+1} \circ \dd^l\!\big)(f)(t) = \sum_{C} \sum_{j=0}^{l}  (-1)^j \sgn(C) \, f(t_{\delta_j(C)}) + \sum_{C}  (-1)^{l+1} \sgn(C) f(t_{C \smallsetminus \{[l]\}})
\end{equation}
holds, where $C$ runs over $\Delta_{l+1}(l+1)$ in both sums. We now consider both sums on the right hand side separately, using \autoref{rem:star}. We claim that the first sum vanishes whereas the second sum equals  $\partial_{\emptyset}f(t) - \sum_{i=0}^l (-1)^i (\partial^i \circ \beta^l)(f)(t)$. 

Concerning the first sum, we recall from the remark that for $C$ and $j \in [l]$ appearing in the sum there exist unique $C' \in \Delta_{l+1}(l+1)$ and $m = m(j)$ such that $\delta_j(C) = \delta_m(C')$. It thus suffices to show that for every such $C$ and $j$,
\begin{equation} \label{signaturee}
	(-1)^j \sgn(C) + (-1)^m \sgn(C') = 0. 
\end{equation}
For this let $C \in \Delta_{l+1}(l+1)$ be as in \eqref{stdform}. According to \autoref{rem:star}, there are three cases we have to consider, depending on the value of $j$: 
\begin{enumerate}[leftmargin=58pt,itemsep=3pt]
	\item[$(j\!<\!k+1) \ $] $(-1)^m \sgn(C') = (-1)^k \cdot (-1)^k \sgn(\phi(C')) = \sgn((s_j \, s_{j+1}\cdots \, s_k)) \cdot  \sgn(\phi(C))$
	
	$= (-1)^{k-j} \cdot (-1)^{k+1} \sgn(C) = (-1)^{j+1} \sgn(C)$, where $(s_j \, s_{j+1}\cdots \, s_k) \in \Sym_{l}$. 
	\item[$(j\!=\!k+1) \ $] $(-1)^m \sgn(C') = (-1)^{k+1} \cdot (-1)^{k+2} \sgn(\phi(C')) = -\sgn(\phi(C)) = (-1)^{j+1} \sgn(C).$
	\item[$(j\!>\!k+1) \ $] $(-1)^m \sgn(C') = (-1)^j \cdot (-1)^{k+1} \sgn(\phi(C')) = (-1)^{j+k+1} \sgn((s_{j-1} \ s_j)) \sgn(\phi(C))$ 
	
	$= (-1)^{j+1} \sgn(C)$, where $(s_{j-1} \ s_j) \in \Sym_{l}$. 
\end{enumerate}
We see that in each of these cases, \eqref{signaturee} holds, and hence the first sum in \eqref{computationhomot} vanishes. Concerning the second sum, we first consider the special chain $C_0 = (0 < \cdots < l < [l])$ and observe that
$	
	(-1)^{l+1} \sgn(C_0) \, f(t_{(0<\cdots<l)}) = (-1)^{l+1} \cdot (-1)^{l+1} \cdot (f \circ \ev_{(0< \cdots < l)})(t) = \partial_\emptyset f(t).
$
It thus remains to show that \vspace{-3pt}
\[
 \sum_{C \neq C_0}  (-1)^{l+1} \sgn(C) f(t_{C \smallsetminus \{[l]\}}) =  - \sum_{i=0}^l (-1)^i (\partial^i \circ \beta^l)(f)(t).
\]
For this we note that if $C \neq C_0$, then there exists some $i \in [l]$ such that $B := C \setminus\{[l]\}  \in \Delta_l([l] \smallsetminus \{i\})$, and one verifies the sign identity
\[
	\sgn(B) = \sgn((i \, \cdots \, l)) \cdot \sgn(C) = (-1)^{l-i} \, \sgn(C),
\]
where $(i \, \cdots \, l) \in \Sym_l$. Thus, 
\begin{align*}
&\sum_{C \neq C_0}  
	(-1)^{l+1} \sgn(C) f(t_{C \smallsetminus \{[l]\}}) = \sum_{i=0}^l \sum_{B \in \Delta_l([l] \smallsetminus \{i\})} \!\!(-1)^{l+1}  (-1)^{l-i} \sgn(B) \ f(t_{B}) \\
	&= - \sum_{i=0}^l (-1)^i \sum_{B \in \Delta_l([l] \smallsetminus \{i\})} \sgn(B) \ f(t_{B}) = - \sum_{i=0}^l (-1)^i \, (\partial^i \circ \beta^l)(f)(t). \qedhere
\end{align*} 
\end{proof}
Using the subdivision maps, we can now establish the main theorem of this section.
\begin{proof}[Proof of \autoref{thm:vanishing}]
For every $l \in [\gamma+1]$, let $\sigma_l \in \Prob(\scD_l)$ as in \autoref{PBSMeasures}. The disintegration $\sigma^{(-)}$, as in \eqref{sigma_disint}, of $\sigma_{l+1}$ along the fibers of the p.m.p map $\partial_\emptyset: \scD_{l+1} \to \BX_l$ produces the bounded operator
$
	\int_{\scD_{l+1}} \dd\!\sigma^{(-)}: \Linfty(\scD_{l+1}) \to \Linfty(\BX_l). 
$
For all $l \in [\gamma+1]$, we set
\begin{equation} \label{defn_homot}
	h^l: \Linfty(\BX_{l+1}) \to \Linfty(\BX_l), \quad h^l\!f(\mathbf{p}):= \int_{\scD_{l+1}} \beta^{l+1}\!f \, \dd\!\sigma^{\mathbf{p}},
\end{equation}
and verify the homotopy relation $h^{l} \circ \dd^{l} + \dd^{l-1} \circ h^{l-1} = \id$. For any $f \in \Linfty(\BX_l)$ and almost all $\mathbf{p} \in \BX_{l}$, we have $(h^{l} \circ \dd^l)(f)(\mathbf{p}) = \int_{\scD_{l+1}} (\beta^{l+1}\circ \dd^l)(f) \, \dd\!\sigma^{\mathbf{p}}$; \vspace{-3pt}
\begin{align}
	f(\mathbf{p}) &= \int_{\BX_l} f \, \dd(\mathrm{Dirac}_{\mathbf{p}}) = \int_{\scD_{l+1}} f \, \dd(\partial_{\emptyset \, \ast} \sigma^{\mathbf{p}}) = \int_{\scD_{l+1}} \partial^\emptyset\! f \, \dd\!\sigma^{\mathbf{p}};\nonumber \\[-4pt] 
	(\dd^{l-1} &\circ h^{l-1})(f)(\mathbf{p}) = \sum_{i=0}^{l} (-1)^i \int_{\scD_l} \beta^l\!f \, \dd\!\sigma^{\delta_i(\mathbf{p})} = \sum_{i=0}^{l} (-1)^i \int_{\scD_{l+1}} \beta^l\!f \, \dd(\partial_{i\, \ast}\sigma^{\mathbf{p}}) \label{eq:dh} \\[-8pt] 
	&= \int_{\scD_{l+1}} \, \sum_{i=0}^{l} (-1)^i \,  (\partial^{i} \circ \beta^l)(f) \,  \dd\!\sigma^{\mathbf{p}}.\nonumber 
\end{align}
The second equality in \eqref{eq:dh} is \autoref{partial_fiber}. Now \autoref{thm:combinatorial} completes the proof. 
\end{proof}
At this point we have established \autoref{thm:vanishing} (and thus \autoref{AcyclicityBound}) assuming the existence of measures $\mu_l$ and $\sigma_l$ as in \autoref{PBSMeasures}. This assumption will be justified in the next section


\section{Measures on configuration spaces}
The purpose of this section is to establish \autoref{PBSMeasures} and thereby to finish the proof of \autoref{AcyclicityBound}. We keep the notation of the previous section and proceed to construct measures $\mu_l$ and $\sigma_l$ as demanded by \autoref{PBSMeasures}.

%

\subsection{Random points on linear sections of $\clP$ and the measures $\mu_l$} \label{subsec:perps}
For a compact, separable metric space $(X,d)$, let $\scC(X)$ be the (compact, separable, metrizable) space of closed subsets of $X$ with the Chabauty topology.  Since $X$ is compact, this topology coincides with the topology on $\scC(X)$ induced by the Hausdorff distance. For any $s \in \bbR_{\geq 0}$, we write $\clH^s$ for the $s$-dimensional Hausdorff measure on $X$ (see \cite[\S 2.10.2]{Federer}), and $\clH^s \mres E$ for its restriction to a Borel subset $E \subset X$, i.e. 
\[
	\clH^s\mres E(A):= \clH^s(A \cap E) \quad \mbox{for any Borel subset } A \subset X. 
\]
Let us write $\dim_\clH$ for the Hausdorff dimension of Borel subsets of $X$; see \cite{Mattila}. We will resort to the following measurability statements. 

\begin{lem} \label{thm:lemma_borel}In the setting above, the following maps are Borel:
\begin{enumerate}[label=\emph{(\roman*)},leftmargin=25pt]
\item $K(-): \scC(X) \to \scC(X), \ A \mapsto K(A) := A \cap K$, for fixed $K \in \scC(X)$.
\item $\clH^s: \scC(X) \to [0,+\infty], \ A \mapsto \clH^s(A)$, for fixed $s \in \bbR_{\geq 0}$. 
\item $\dim_\clH: \scC(X) \to [0,+\infty], \ A \mapsto \dim_\clH(A)$.
\end{enumerate}
\end{lem} 

\begin{proof}
The functions from (ii) and (iii) are of Baire class 2, hence Borel; see \cite{Mattila}. On the other hand, one verifies easily that the function (i) is \emph{upper semicontinuous} in the following sense: For any sequence $(A_n)_n \subset \scC(X)$ with limit $A \in \scC(X)$, the inclusion $\limsup_n (A_n \cap K) \subset A \cap K$ holds, where $\limsup_n (A_n \cap K) \in \scC(X)$ is the set of cluster points of sequences $(a_n) \in X$ with $a_n \in A_n \cap K$. It is well known that maps with this property are Borel; see e.g. \cite{LeBoudec-MatteBon}. 
\end{proof}

Let us now fix an auxiliary inner product on $V$. Its pullback to the unit sphere $\bbS(V)$ descends ---after realification if $\kay = \bbC$--- to a Riemannian metric on $\bbP(V) \cong \bbS(V)/\bbS(\kay)$. Since $\GL(V)$ acts transitively on each Grassmannian $\Gr_k(V)$, it follows from the open mapping theorem for homogeneous spaces that the subspace topology on $\Gr_k(V) \subset \scC(\bbP(V))$ coincides with the quotient topology with respect to $\GL(V)$; in the sequel we will always equip $\Gr_k(V)$ with this topology and $\Gr(V)$ with the corresponding direct sum topology. 

We recall that $\clP = \BX_0 = X_0 \subset \bbP(V)$ denotes the set of isotropic points in $\bbP(V)$, which is a real projective algebraic set (and coincides with $\bbP(V)$ if and only if $\omega$ is symplectic). Now let $W \in \Gr(V)$ so that $W \subset \bbP(V)$. We denote by $\clP(W)$ the real projective algebraic set $\clP(W) := W \cap \clP\subset \bbP(V)$. We then set
\[
	\Gr_k^\omega(V) := \{W \in \Gr_k(V) \ \mid \ \clP(W) \neq \emptyset\}, \quad\mbox{resp.} \quad \Gr^\omega(V) := {\textstyle \bigsqcup_k} \Gr_k^\omega(V).
\]
Observe that these are closed, $G$-invariant subsets of $\Gr_k(V)$ resp. $\Gr(V)$. For $W \in \Gr^\omega(V)$, let $\clP(W)_{\max}$ be the smooth locus of the union of all Zariski-irreducible components of $\clP(W)$ of dimension $d(W):= \dim \clP(W)$. In other words, $\clP(W)_{\max}$ is obtained from $\clP(W)$ by removing all components of non-maximal dimension and all singular points from the components of maximal dimension. This ensures that $\clP(W)_{\max}$ is an embedded submanifold of $\bbP(V)$ and that its induced Riemannian measure corresponds (up to rescaling) to the $d(W)$-dimensional Hausdorff measure; see \cite[\S 3.2.46]{Federer}. Thus, $\clH^{d(W)}(\clP(W)) = \clH^{d(W)}(\clP(W)_{\max}) > 0$, and $\clH^{d(W)}(\clP(W)) < \infty$ by compactness. In particular, $\dim_\clH(\clP(W)) = d(W) \in \bbN$. 

\begin{defn}
For every $W \in \Gr^\omega(V)$, we define the probability measure $\lambda_W$ on $\bbP(V)$ as
\[
	\lambda_W := \frac{\clH^{d(W)} \mres \clP(W)}{\clH^{d(W)}(\clP(W))}. 
\]
\end{defn}

\begin{lem} \label{thm:borel_Gqi} 
The map $\lambda: \Gr^\omega(V) \to \Prob(\bbP(V))$ defined by the expression above is \emph{(i)} such that $\lambda_W(\clP(W))=1$ for all $W$, \emph{(ii)} Borel, and \emph{(iii)} \emph{$G$-quasi-equivariant}, i.e. the measures $g_\ast \lambda_W$ and $\lambda_{gW}$ are mutually absolutely continuous for every $g \in G$ and $W \in \Gr^\omega(V)$. 
\end{lem}

\begin{proof}
Point (i) is obvious. For (ii), it is enough to check that $h_K: W \mapsto \clH^{d(W)}(A \cap \clP(W))$ is Borel for any closed $K \subset \bbP(V)$. This is true since
$
	h_K = \clH \, \circ \, (\dim_\clH \, \times \, K(-)) \, \circ \, \clP(-) \, \circ \, \iota,
$ with
$\clH: \bbN \ \times \ \scC(\bbP(V)) \to [0,+\infty]$ being the map $(d,C) \mapsto \clH^d(C)$;
$\dim_\clH\!\!: \scC(\bbP(V)) \to [0,+\infty]$, $K(-), \clP(-): \scC(\bbP(V)) \to \scC(\bbP(V))$ being as in \autoref{thm:lemma_borel}; and
$\iota: \Gr^\omega(V) \hookrightarrow \scC(\bbP(V))$ being the inclusion. 
These are Borel by \autoref{thm:lemma_borel}. The last step in the composition is well defined, since the image of $\dim_\clH \, \circ \, \clP(-) \, \circ \, \iota$ lies in $\bbN$. 

For (iii), it suffices to show that $g_\ast (\clH^{d(W)}\mres \clP(W))$ and $\clH^{d(gW)} \mres \clP(gW)$ are mutually absolutely continuous for every $g \in G$ and $W \in \Gr_k^\omega(V)$. Indeed, since $g: \clP(W) \to \clP(gW)$ is bi-Lipschitz, then $d(W) = d(gW)$, and by the area formula \cite[\S 3.2.3, \S 3.2.46]{Federer},
\[
	g_\ast (\clH^{d(W)}\mres \clP(W))(A) = \int_{\clP(W)} \bbone_{A}(gx) \, \dd\!\clH^{d(W)}(x) = \int_{\clP(gW)} \bbone_{A}(x) \cdot (\clJ g^{-1})(x)\, \dd\!\clH^{d(gW)}(x); 
\]
here $\clJ g^{-1}$ denotes the Jacobian of $g^{-1}$, computed with respect to the induced Riemannian metrics on $\clP(W)$ and $\clP(gW)$. The result then follows since this Jacobian is invertible.
\end{proof}

We record the following property of the measures $\lambda_W$ for later use.

\begin{lem} \label{lemmaaaaa} Assume that $r \geq 3$. For $i \in \{1,2\}$, let $A_i \in \Gr^\omega(V)$ be either a maximal isotropic subspace, or a subspace of codimension at most $r-2$. In addition, assume that $A_1 \neq A_2$ if both $A_1$ and $A_2$ are maximal isotropic. Then $\omega(a_1,a_2) \neq 0$ for $(\lambda_{A_1} \otimes \lambda_{A_2})$-almost all $(a_1, a_2) \in A_1 \times A_2$.
\end{lem}

\begin{proof} Let $W_i < V$ be a linear subspace with $A_i = \bbP(W_i)$ for $i = 1,2$. We claim that $\Span(W_i \cap \clV) = W_i$ and that $\clP(A_i)$ is a projective variety. This is immediate if $A_i$ is totally isotropic, since $\clP(A_i) = A_i$ and $W_i \cap \clV = W_i$; if instead the codimension hypothesis holds, the claim follows from \autoref{linalg3}. Our second claim is that the set
\[
M := \{(v_1, v_2) \in (W_1 \cap \clV) \times (W_2 \cap \clV) \mid \omega(v_1, v_2) \neq 0 \}
\]
is non-empty. If $W_1$ and $W_2$ are two distinct maximal isotropic subspaces, this is again immediate. Hence, assume that $A_1$ has codimension at most $r-2$. We have
\[
	\dim W_1 \geq r+d+2 \qand \dim W_2 \geq r.
\]
In particular, $\dim W_1 + \dim W_2 > 2r + d = \dim V$. If the set $M$ were empty, then for all $v_1 \in W_1 \cap \clV$, we would have $W_2 \cap \clV \subset v_1^\perp$ and thus $W_2 = \Span(W_2 \cap \clV) \subset v_1^\perp$, and hence $v_1 \in W_2^\perp$. This would imply that $W_1 \cap \clV \subset W_2^\perp$ and hence $W_1 = \Span(W_1 \cap \clV) \subset W_2^\perp$. In terms of dimension, this would mean that $\dim W_1 + \dim W_2 \leq \dim V$, a contradiction. 

We have shown that 
$
\{(p_1, p_2) \in \clP(A_1) \times \clP(A_2) \mid \omega(p_1, p_2) \neq 0\}
$
is a non-empty Zariski-open subset in the real projective variety $\clP(A_1) \times \clP(A_2)$. Because of the irreducibility of the $\clP(A_i)$'s, the measure $\lambda_{A_1} \otimes \lambda_{A_2}$ belongs to the Lebesgue class of $\clP(A_1) \times \clP(A_2)$, so the complement of the subset above is null. 
\end{proof}

As a first application of the assignment $\lambda$ we can now construct the measures $\mu_l$ from \autoref{PBSMeasures}. Recall from Subsection \ref{subsec:grass} that for $k \in [r-1]$ we denote by $\clG_k \subset \Gr_k(V)$ the subset formed by all totally isotropic projective subspaces of dimension $k$. Our construction depends on the choice of a $G$-quasi-invariant probability measure $\nu \in \Prob(\clG_{r-1})$. There is complete freedom in the choice of this measure (and hence the measures $\mu_l$ will be far from unique), but it is necessary to fix a choice of $\nu$ once and for all.

\begin{con}\label{ConMu} Given $l \in [r]$, let $\clF_l$ be the set of tuples $(p_0,\ldots,p_l; W) \in \bbP(V)^{l+1} \times \clG_{r-1}$ such that $p_0,\ldots,p_l \in W$ and denote the canonical projections by
$\mathrm{proj}_{ \bbP(V)^{l+1}}: \clF_l \to   \bbP(V)^{l+1}$, $(p_0,\ldots,p_l; W) \mapsto (p_0, \ldots, p_l)$ and $\mathrm{proj}_{\clG_{r-1}}: \clF_l \to  \clG_{r-1}$, $(p_0,\ldots,p_l; W) \mapsto W$ respectively. The fiber of $\mathrm{proj}_{\clG_{r-1}}$ over $W \in \clG_{r-1}$ is given by $W^{l+1}\times \{W\}$. We now endow this fiber with the product measure $\lambda_W^{\otimes l+1} \otimes \mathrm{Dirac_W}$, and integrate to obtain a measure
\[
	\tilde{\mu}_l := \int_{\clG_{r-1}} (\lambda_W^{\otimes l+1} \otimes \mathrm{Dirac}_W) \, \dd\!\nu(W) \ \in \ \Prob(\clF_l);
\]
we then define  $\mu_l := (\mathrm{proj}_{ \bbP(V)^{l+1}})_\ast \tilde{\mu}_l \in \Prob(\BX_l)$. 
\end{con}
\begin{lem} The measures $\mu_l \in \Prob(\BX_l)$ are $G$-quasi-invariant and satisfy Properties (i) and (ii) from \autoref{PBSMeasures}.
\end{lem}
\begin{proof} The measures $\tilde{\mu}_l  \in  \Prob(\clF_l)$ are well-defined and $G$-quasi-invariant by \autoref{thm:borel_Gqi} and quasi-invariance of $\nu$, and since the map $\mathrm{proj}_{ \bbP(V)^{l+1}}$ is $G$-equivariant and takes values in $\BX_l$ the measures $\mu_l$ are well-defined and $G$-quasi-invariant. Now Properties (i) and (ii) follows from Fubini's theorem, taking into account the product structure of the measure $\lambda_W^{\otimes l+1}$.
\end{proof}

A second application of $\lambda$ will be the formalization of the auxiliary measures $\lambda_\perp^{(-)}$, approached heuristically in Subsection \ref{heuristic}. Consider the maps
\[\begin{array}{ll}
     (-)^\perp\, : \  \Gr(V) \to \Gr(V),  & (\bbP(L))^\perp := \bbP(L^\perp), \qand \\[3pt]
     \Span: \ \clP^{k+1} \to \Gr(V) & \Span(p_0,\ldots,p_k) := \Span\{p_0,\ldots,p_k\}, 
\end{array}\]
for any $k \in [n-2]$. Both are $G$-equivariant and Borel, the former being continuous and the latter piecewise continuous. Recall that $(-)^\perp$ is an involution. 

\begin{defn} \label{def_randomperp}
For any $k \in \bbN_{>0}$, we define $\Omega_k := \{\mathbf{p} \in \clP^{k+1} \, \mid \, \Span(\mathbf{p})^\perp \in \Gr^\omega(V)\}$, and  
\[
	\lambda_\perp = \lambda_{\perp,k}: \Omega_k \to \Prob(\bbP(V)), \quad \lambda_\perp^{\mathbf{p}} := \lambda_{\Span(\mathbf{p})^\perp}.
\]
Note that $\Omega_k \neq \emptyset$ is $G$-invariant and Borel. Moreover, by \autoref{thm:borel_Gqi}, the map $\lambda_\perp$ is Borel, $G$-quasi-equivariant, $\Sym_k$-invariant, and such that $\lambda_\perp^{\mathbf{p}}(\clP(\Span(\mathbf{p})^\perp)) = 1$ for every $\mathbf{p} \in \Omega_k$. 
\end{defn}

\begin{rem} \label{everythingfits}
Observe that if $k + 1< 2r$, then $\Omega_k = \clP^{k+1}$. Indeed, if $\mathbf{p} \in \clP^{k+1}$, then the subspace $\Span(\mathbf{p})^\perp$ has codimension strictly less than $2r$, hence the underlying subspace of $V$ has dimension greater than $d$. By \autoref{thm:adapted_basis}, $\Span(\mathbf{p})^\perp$ must contain an isotropic point. We thus obtain a well-defined map $\lambda_\perp: \clP^{k+1} \to  \Prob(\bbP(V))$.

More generally, let $S$ be a set of cardinality $|S| = k+1 < 2r$. If we fix an enumeration $S = \{s_0, \dots, s_{k}\}$, then we may define
\[
\lambda_\perp: \clP^{S} \to \Prob(\bbP(V)), \quad \lambda_\perp^{\mathbf p} := \lambda_\perp^{(\mathbf p(s_0), \dots, \mathbf p(s_k))}.
\]
This map is actually independent of the choice of enumeration by the aforementioned invariance under permutations.
\end{rem}

\subsection{Construction of measures on $\scD_l$}\label{SecMeasuresPBS} 
We work in the same setting as in the previous subsection; in particular, $\nu \in \Prob(\clG_{r-1})$ is a fixed $G$-quasi-invariant probability measure and the measures $\mu_0, \dots, \mu_r$ are constructed from $\nu$ as in \autoref{ConMu}. The goal of this subsection is to complete the proof of \autoref{PBSMeasures} and in particular to construct the measures $\sigma_l$. The nature of our construction will be inductive. For its purpose, we introduce a filtration on the vertex set $\Delta_0 S$ of the complex $\Delta S$ for any finite totally ordered set $S$. 
\begin{con}\label{ProjectionsEpsiloni} Let $S$ be a finite totally ordered set. For any $c \in \Delta_0 S = S \cup 2^S$, we let
\[
|c| := \left\{\begin{array}{rl} 0, & \text{if } c\in S\\ 
\mathrm{card}(c), & \text{if } c \in 2^S
 \end{array}\right..
\]
We then set
\begin{equation*}
	\Delta_0S_m := \{c \in \Delta_0 S \, \mid \, |c| \leq m\} \quad \mbox{for } m \in [|S|].
\end{equation*}
Clearly $\Delta_0S_0 \subset \Delta_0S_1 \subset \cdots \subset \Delta_0S_{|S|} = \Delta_0 S$, and the difference $\partial\Delta_0S_{m+1} := \Delta_0S_{m+1}  \smallsetminus\Delta_0S_{m} $ consists of all subsets of $S$ of cardinality $m+1$. Moreover, these filtrations are functorial in the sense that if $T \subset S$, then there are obvious inclusions $\Delta_0T_{m}  \hookrightarrow \Delta_0S_{m} $ for every $m \in [|T|]$. We also observe that the inclusions $\Delta_0S_{m} \hookrightarrow\Delta_0S_{m+1}$ and $\Delta_0T_m \hookrightarrow \Delta_0S_m$ induce $G$-equivariant projections
\begin{equation} \label{proj_eta} 
	\pi_m: \ \clP^{\Delta_0S_{m+1}} \to \clP^{\Delta_0S_m} \qand \eta_{T\subset S}: \ \clP^{\Delta_0S_m}\to \clP^{\Delta_0T_m}   
\end{equation}
respectively. As before we are particularly interested in the case where $S = [l-1]$ for some $l \in \bbN$. In this case, we write $\Delta_0(l)_{m} := \Delta_0S_m$ so that
\[
\Delta_0(l)_0 \subset \Delta_0(l)_1 \subset \dots \subset \Delta_{0}(l)_l = \Delta_0(l).
\]
If $i \in [l-1]$, then there is a unique order-preserving bijection $[l-1] \to [l] \smallsetminus \{i\}$, and we denote by $\varepsilon_i: \Delta_0(l)_m \to \Delta_0([l] \smallsetminus \{i\})_m$ be the induced $G$-equivariant homeomorphism on the level of filtrations.
\end{con}
Using the above filtration we will construct the desired measures $\sigma_l \in \Prob(\clD_l)$ in three steps:
\begin{enumerate}
\item We first construct for any finite totally ordered set $S$ a family of $G$-quasi-invariant probability measures $\sigma_{S,m} \in \Prob(\clP^{\Delta_0S_m})$ by induction on $m$. In particular, for $S = [l-1]$, this will yield a family of measures $\sigma_{l,m} := \sigma_{[l-1], m} \in \Prob(\clP^{\Delta_0(l)_m})$. 
\item We then show (again by induction) that the measure $\sigma_{l,l} \in  \Prob(\clP^{\Delta_0(l)})$ gives full measure to the subspace $\scS(l)$ from \eqref{scS}.
\item Finally, we use the identification $\mathrm{res}: \scD_l \to \scS(l)$ from \autoref{rem:incarnations} to define $\sigma_l := \mathrm{res}^{-1}_*\sigma_{l,l} \in \Prob(\scD_l)$.
\end{enumerate}
From now on, we work in the setting of \autoref{PBSMeasures}; in particular, $l \in [\gamma+1]$. Moreover,  $S = \{s_0 < \cdots < s_{l-1}\}$ is a totally ordered set of cardinality $l$. We construct $\sigma_{S,0}$:
\begin{con}\label{Con1a} On $\clP^{l+1} = \clP^l \times \clP^1$ we consider the product measure $\mu_{l-1} \otimes \mu_0$. We then define $\sigma_{S, 0}$ as the pushforward measure of this measure under the $G$-equivariant homeomorphism $\clP^{l+1} \to \clP^{\Delta_0S_0}$ which identifies $(t_0, \dots, t_{l})$ and $(t_{s_0}, \dots, t_{s_{l-1}}, t_\emptyset)$, i.e.
\[
\dd\!\sigma_{S,0}(t_{s}\mid s \in \Delta_0S_0) = \dd\!\mu_{l-1}(t_{s_0}, \dots, t_{s_{l-1}})\; \dd\!\mu_0(t_\emptyset).
\]
Note that the measure $\sigma_{S,0} \in \Prob(\clP^{\Delta_0S_0})$ is $G$-quasi-invariant.
\end{con}
We construct the higher measures $\sigma_{S,m}$ inductively. For this we assume that $G$-quasi-invariant measures $\sigma_{S, 0}, \dots, \sigma_{S,m}$ have been defined for some $m \in [l-1]$. We then define a map  \[\tilde\sigma_{S,m+1}: \clP^{\Delta_0S_{m}} \to \Prob(\clP^{\Delta_0S_{m+1}})\] as follows: For any $\varphi \in C(\clP^{\Delta_0S_{m+1}})$ and $(t_c \, \mid \, c \in \Delta_0S_{m}) \in \clP^{\Delta_0S_{m}}$, we set
	 \begin{equation} \label{definduction}
	 	\tilde\sigma_{S,m+1}(t_c \!\mid\! c \in \Delta_0S_{m})(\varphi) \!:= \!\int\! \varphi(t_c \!\mid\! c \in \Delta_0S_{m+1}) \, \dd\!\Big(\!\bigotimes_{c_0} \lambda_\perp^{(t_c \mid c < c_0)}\!\Big)\big(t_{c_0} \!\mid\!  c_0 \in \partial\Delta_0S_{m+1}\big).
		\end{equation}
	\begin{lem} \label{thm_welldefined}
	The formula \eqref{definduction} gives a well-defined map $\tilde\sigma_{S,m+1}: \clP^{\Delta_0S_{m}} \to \Prob(\clP^{\Delta_0S_{m+1}})$ that is Borel and $G$-quasi-equivariant. 
	\end{lem}
	\begin{proof} To see that $\tilde\sigma_{S,m+1}$ is well-defined, we have to check that $\lambda_\perp^{(t_c \mid c < c_0)}$ is well-defined for all $c_0 \in  \partial\Delta_0S_{m+1}$. In view of \autoref{everythingfits}, this amount to showing that $|\{c \mid c < c_0\}| < 2r$ for all $c_0 \in  \partial\Delta_0S_{m+1}$. If $\tilde\gamma_\ast:\bbN \to \bbN$ is the function from \eqref{defgamma}, then we have for all such $c_0$ that 
\begin{equation} \label{annoying_bound}
\begin{array}{rcl}
|\{c \mid c < c_0\} | &=& |\Delta_0(c_0) \smallsetminus \{c_0\}| = m+2^{m+1} = 2\left(\frac{m+1}{2} + 2^m\right) - 1 \leq 2\tilde\gamma_\ast(m) - 1 \\
& \leq & 2\tilde\gamma_\ast(l-1) - 1 \leq r-2
\end{array}
\end{equation}
where the last inequality is equivalent to the assumption from \autoref{PBSMeasures} that $l - 1 \leq \tilde\gamma(\lfloor(r-1)/2\rfloor) = \gamma$. In particular, $|\{c \mid c < c_0\}| < 2r$, and hence, $\tilde\sigma_{S,m+1}$ is well-defined.
		
To see that $\tilde\sigma_{S,m+1}$ is Borel and $G$-quasi-equivariant, we observe that for every $(t_c \mid c \in \Delta_0S_{m})$, the evaluation $\tilde\sigma_{S,m}(t_c \mid c \in \Delta_0S_{m})(f)$ is integration of $f \in C(\clP^{\Delta_0S_{m+1}})$ against 
	\[
		\mathrm{Dirac}_{(t_c \mid c \in \Delta_0S_{m})} \, \otimes \, \bigotimes_{c_0 \in \partial \Delta_0S_{m+1}} \lambda_\perp^{\eta_{c_0 \, \subset \, S}(t_c \mid c \in\Delta_0S_{m})}
	\]
	under the $G$-equivariant identification $\clP^{\Delta_0S_{m+1}} \cong \clP^{\Delta_0S_{m}} \times \clP^{\partial \Delta_0S_{m+1}}$.
	\end{proof}
\begin{con}\label{Con1b} Using the function $\tilde\sigma_{S,m+1}: \clP^{\Delta_0S_{m}} \to \Prob(\clP^{\Delta_0S_{m+1}})$ we now define
\[
		\sigma_{S,m+1} := \int \tilde\sigma_{S,m+1} \, \dd\!\sigma_{S,m} \in \Prob(\clP^{\Delta_0S_{m+1}}),\]
where the integral is in the sense of \eqref{defdisint}. Note that the measure $\sigma_{S,m+1} \in \Prob(\clP^{\Delta_0S_{m+1}})$ is $G$-quasi-invariant by \autoref{thm_welldefined} and the induction hypothesis.
\end{con}
At this point we have completed Step (1) of our construction. Recall that for $S := [l-1]$ we write $\sigma_{l,m} := \sigma_{[l-1],m}$. The next example should clarify the nature of these measures.

\begin{exmpl} \label{exmpl_measures}
Let $l \in \{1,2\}$, assume $r \geq \gamma_\ast(l-1)$, and let $S = [l-1]$. If $l = 1$, then $r \geq 2$, $\Delta_0(1)_{0} = \{0,\emptyset\}$ and $\partial \Delta_0(1)_{1} = \{\{0\}\}$. The measures from \autoref{Con1a} and \autoref{Con1b} are $\sigma_{1,0} = \mu_0 \otimes \mu_0$, and 
\[
	\sigma_{1,1}(\varphi) = \iint \left(\int \varphi(t_0,t_\emptyset,t_{\{0\}})\, \dd\!\lambda_\perp^{t_0,t_\emptyset}(t_{\{0\}})\right) \dd\!\mu_0(t_0) \, \dd\!\mu_0(t_\emptyset).
\]
Analogously, if $l = 2$, then $r \geq 3$, $\Delta(2)_0 = \{0,1,\emptyset\}$, $\partial \Delta_0(2)_1 = \{\{0\}, \{1\}\}$, and $\partial \Delta_0(2)_2 = \{\{0,1\}\}$. Furthermore, $\sigma_{2,0} = \mu_1 \otimes \mu_0$, 
\begin{equation*}
	\sigma_{2,1}(\varphi_1) = \iint \left(\iint \varphi_1(t_0,t_1,t_\emptyset,t_{\{0\}},t_{\{1\}}) \, \dd\!\lambda_\perp^{t_0,t_{\emptyset}}(t_{\{0\}}) \, \dd\!\lambda_\perp^{t_1,t_{\emptyset}}(t_{\{1\}}) \right) \dd\!\mu_1(t_0,t_1) \, \dd\!\mu_0(t_\emptyset), \mbox{ and}
\end{equation*}
\begin{equation*}
	\sigma_{2,2}(\varphi) = \int \left(\int \varphi(t_0,t_1,t_\emptyset,t_{\{0\}},t_{\{1\}},t_{\{0,1\}}) \, \dd\!\lambda_\perp^{t_0,t_1,t_\emptyset,t_{\{0\}},t_{\{1\}}}(t_{\{0,1\}})\right) \dd\!\sigma_{2,1}(t_0,t_1,t_\emptyset,t_{\{0\}},t_{\{1\}}).
\end{equation*}
\end{exmpl}
Observe that in the last example the points $(t_0,t_\emptyset,t_{\{0\}})$ distributed according to $\sigma_{1,1}$ belong almost surely to the set \eqref{points_delta1} in \autoref{exmpl_delta}, and therefore, can be used to define a probability measure on $\scD_1$. This observation will be generalized in the next step.

For Step (2) of our construction we need to define a filtration of the space $\scS(l)$ which mirrors our filtration of $\Delta_0(l)$. 

\begin{defn} We define $\scS(l)_{-1} := \emptyset$ and inductively
\begin{eqnarray*}
\scS(l)_{0} &:=&\left\{(t_c \mid c \in \Delta_0(l)_{0}) \mid 
			\forall c_0 \in \partial \Delta_0(l)_0 \ \forall (c < c_0) \in \Delta_1 (l) \,: \, (t_c,t_{c_0}) \in \BX_1 
			\right\},\\
\scS(l)_{m+1} &:=& \left\{(t_c \mid c \in   \Delta_0(l)_{m+1}) \left| \begin{array}{l}
			(t_c \mid c \in  \Delta_0(l)_{m}) \in \scS_{m}; \qand\\
			\forall c_0 \in \partial  \Delta_0(l)_{m+1} \ \forall (c < c_0) \in \Delta_1 (l) \,: \, (t_c,t_{c_0}) \in \BX_1 
			\end{array}\right.\right\}.
\end{eqnarray*}
\end{defn}
\begin{lem} \label{lem_support}
For all $l \in [\gamma+1]$ and $m \in [l]$, the set $\scS(l)_m \subset \clP^{\Delta_0(l)_m}$ is closed, and satisfies $\sigma_{l,m}(\scS(l)_m) = 1$.
\end{lem}
\begin{proof}
The second line of conditions in the definition of $\scS(l)_{m+1}$ is a collection of closed conditions, hence the assertion that $\scS(l)_{m}$ is closed for every $m \in [l]$ follows by induction over $m$. Concerning the second claim, we first observe that $\sigma_{S,0}(\scS(l)_0) = 1$ by definition. Now, assume $\sigma_{S,m}(\scS(l)_m) = 1$ has been established for some $m \in [l]$. Then, for any $(t_c \mid c \in \Delta_0(l)_{m+1})$, we have the equality \vspace{-3pt}
\[
	\bbone_{\scS_{m+1}}(t_c \mid c \in \Delta_0(l)_{m+1})= \bbone_{\scS_{m}}(t_c \mid c \in\Delta_0(l)_{m}) \cdot \prod_{c_0 \in \partial \Delta_0(l)_{m+1}} \bbone_{\clP(\Span(t_c \mid c < c_0)^\perp)} (t_{c_0}), 
\]
of characteristic functions, and hence, $\sigma_{l,m+1}(\scS(l)_{m+1}) = \sigma_{l,m}(\scS(l)_{m}) = 1$. 
\end{proof}
\begin{cor} For every $l \in [\gamma+1]$, we have $\sigma_{l,l}(\scS(l)) = 1$.\qed
\end{cor}
\begin{con}\label{SigmaFinal}
For every $l \in [\gamma+1]$ we may now use the $G$-equivariant identification $\mathrm{res}: \scD_l \to \scS(l)$ from \autoref{rem:incarnations} to define a $G$-quasi-invariant measure 
\begin{equation}
\sigma_l := \mathrm{res}^{-1}_*\sigma_{l,l} \in \Prob(\scD_l).
\end{equation}
\end{con}

\begin{exmpl}
Following the definition of the homotopy \eqref{defn_homot}, and combining Examples \ref{exmpl_beta} and \ref{exmpl_measures}, we obtain for $h^0$ the exact expression \eqref{eq:heuristics2} from the heuristic discussion. 
\end{exmpl}

In the remainder of this section we are going to show that the measures $\sigma_l$ satisfy the properties required by \autoref{PBSMeasures}. We start with some general lemmas; throughout  $S = \{s_0 < \cdots < s_{l-1}\}$ is a totally ordered set of cardinality $l \in [\gamma+1]$.

\begin{lem} \label{brute}
For $m \in [l-1]$ and $T \subset S$ the maps \vspace{-3pt}
\[
	\pi_m: (\clP^{\Delta_0S_{m+1}}, \sigma_{S,m+1}) \to (\clP^{\Delta_0S_m},  \sigma_{S,m}) \qand \eta_{T\subset S}: \ (\clP^{\Delta_0S_m}, \sigma_{S,m})\to (\clP^{\Delta_0T_m}, \sigma_{T,m})
\]
from \eqref{proj_eta} are p.m.p. 
\end{lem}
\begin{proof} Let $m \in [l-1]$. For any $\varphi \in C(\clP^{\Delta_0{S}_m})$ and $(t_c \!\mid\! c \in \Delta_0S_m) \in \clP^{\Delta_0S_m}$ we have
\begin{eqnarray*}
&& \pi_{m\, \ast}\big(\tilde\sigma_{S,m+1}(t_c \!\mid\! c \in \Delta_0S_m)\big)(\varphi)\\
&=& \!\int\! \varphi(t_c \!\mid\! c \in \Delta_0S_{m}) \, \dd\!\Big(\!\bigotimes_{c_0} \lambda_\perp^{(t_c \mid c < c_0)}\!\Big)\big(t_{c_0} \!\mid\!  c_0 \in \partial\Delta_0S_{m+1}\big)\\
&=&  \varphi(t_c \!\mid\! c \in \Delta_0S_{m}).
\end{eqnarray*}
This then implies the desired statement about the maps $\pi_m$, since
\begin{eqnarray*}
	(\pi_{m\, \ast} \sigma_{S,m+1})(\varphi) &=& \int \pi_{m\, \ast} \big(\tilde\sigma_{S,m+1}(t_c \!\mid\! c \in  \Delta_0S_m)\big)(\varphi) \, \dd\!\sigma_{S,m}(t_c \!\mid\! c \in \Delta_0S_{m})\\
	&=&  \int  \varphi(t_c \!\mid\! c \in \Delta_0S_{m})  \, \dd\!\sigma_{S,m}(t_c \!\mid\! c \in \Delta_0S_m)  \ = \ \sigma_{S,m}(\varphi).
\end{eqnarray*}
Now let $s \in S$ and let $\eta_s := \eta_{S \smallsetminus \{s\} \subset S}$. To establish the second claim it suffices to show that $\eta_s: (\clP^{\Delta_0S_m}, \sigma_{S,m}) \to (\clP^{\Delta_0(S\smallsetminus \{s\})_m}, \sigma_{S \smallsetminus\{s\}, m})$ is p.m.p. for all  $m \in [l]$, and we will show this by induction on $m$. For $m=0$, by \autoref{Con1a}, we have a commutative diagram
\begin{equation*} 
		\begin{gathered}	
		\xymatrixcolsep{3pc}
		\xymatrix@R=1.5pc{(\clP^{\Delta_0S_0}, \sigma_{S,0}) \ar@{<->}[d]_{\cong} \ar[r]^{\eta_{s}} & (\clP^{\Delta_0(S\smallsetminus\{s\})_0}, \sigma_{S\smallsetminus\{s\}, 0})\ar@{<->}[d]^{\cong} \\
		(\BX_{l-1} \times \clP, \mu_{l-1} \otimes \mu_0)  \ar[r]^{\delta_s \, \times \, \id} & (\BX_{l-2} \times \clP, \mu_{l-2} \otimes \mu_0),}
		\end{gathered}
\end{equation*}
in which the vertical arrows are p.m.p. The claim then follows from the fact that the face maps $\delta_s: (\BX_{l-1}, \mu_{l-1}) \to (\BX_{l-2}, \mu_{l-2})$ are p.m.p. by \autoref{PBSMeasures} (ii).

For the induction step, assume that $\eta_{s\,\ast}\sigma_{S,m-1} = \sigma_{S\smallsetminus\{s\},m-1}$ for some $m \in [l-1]$. For any $\varphi \in C(\clP^{\Delta_0(S \smallsetminus\{s\})_m})$ we then have, on the one hand,
\[
\eta_{s\, \ast} \sigma_{S,m}(\varphi) \; =\;\int \tilde \sigma_{S,m}(t_c\mid c\in \Delta_0(S)_{m-1})(\eta_s^*\varphi)\, \dd \sigma_{S, m-1}(t_c\mid c\in \Delta_0S_{m-1})
\]
and on the other hand
\begin{eqnarray*}
\sigma_{S \setminus\{s\}, m}(\varphi) &=& \int \tilde\sigma_{S \setminus\{s\}, m}(t_c \mid c \in \Delta_0 (S \smallsetminus\{s\})_{m-1})(\varphi)\, \dd\sigma_{S \smallsetminus\{s\}, m-1}(t_c \mid c \in \Delta_0 (S \smallsetminus\{s\})_{m-1})\\
&=& \int \tilde\sigma_{S \setminus\{s\}, m}(t_c \mid c \in \Delta_0 (S \smallsetminus\{s\})_{m-1})(\varphi)\, \dd \eta_{s\,\ast}\sigma_{S, m-1}(t_c \mid c \in \Delta_0 (S \smallsetminus\{s\})_{m-1})\\
&=&  \int \tilde\sigma_{S \setminus\{s\}, m}(t_c \mid c \in \Delta_0 (S \smallsetminus\{s\})_{m-1})(\varphi)\, \dd\sigma_{S, m-1}(t_c \mid c \in \Delta_0 S_{m-1}),
\end{eqnarray*}
hence we need to show that almost surely
\[
  \tilde\sigma_{S,m}(t_c\mid c\in \Delta_0(S)_{m-1})(\eta_s^*\varphi) \;  = \;\tilde\sigma_{S \setminus\{s\}, m}(t_c \mid c \in \Delta_0 (S \smallsetminus\{s\})_{m-1})(\varphi).
\]
This follows from
\begin{align*}
	&  \tilde\sigma_{S,m}(t_c\mid c\in \Delta_0(S)_{m-1})(\eta_s^*\varphi) \; = \; \int \eta_s^*\varphi(t_c \!\mid\! c \in \Delta_0S_{m}) \, \dd\tilde\sigma_{S,m}(t_c \!\mid\! c \in \Delta_0S_{m})\\
	=& \int (\eta_s^*\varphi)(t_c \mid c \in \Delta_0S_m) \dd\!\Big(\!\bigotimes_{c_0} \lambda_\perp^{(t_c \mid c < c_0)}\!\Big)\big(t_{c_0} \!\mid\!  c_0 \in \partial\Delta_0S_{m}\big) \\
	=& \!\int\! \varphi(t_c \!\mid\! c \in \Delta_0(S\smallsetminus\{s\})_{m}) \, \dd\!\Big(\!\bigotimes_{c_0} \lambda_\perp^{(t_c \mid c < c_0)}\!\Big)\big(t_{c_0} \!\mid\!  c_0 \in \partial\Delta_0S_{m}\big) \\
	 =& \ \!\int\! \varphi(t_c\mid c \in \Delta_0(S\smallsetminus\{s\})_{m}) \, \dd\!\Big(\!\bigotimes_{c_0} \lambda_\perp^{(t_c \mid c < c_0)}\!\Big) \big(t_{c_0} \!\mid\!  c_0 \in \partial\Delta_0(S\smallsetminus\{s\})_{m}\big) \\
	 =& \int \varphi(t_c \!\mid\! c \in \Delta_0(S\smallsetminus\{s\})_{m}) \, \dd\tilde\sigma_{S\smallsetminus\{s\},m}(t_c \!\mid\! c \in \Delta_0(S\smallsetminus\{s\})_{m})\\
	 	 =& \ \tilde\sigma_{S \setminus\{s\}, m}(t_c \mid c \in \Delta_0 (S \smallsetminus\{s\})_{m-1})(\varphi).\qedhere
\end{align*} 
\end{proof}

\begin{lem} \label{interseccion}
For $l \in [\gamma+1]$, let $C \in \Delta_l(l)$ and $B:= \delta_l(C) \in \Delta_{l-1}(l)$. Then
\[
 	\int\bbone_{\clP(\Span(t_B))}(t_{[l-1]}) \, \dd\!\sigma_l(t) = 0.
\]
\end{lem}
\begin{proof} By definition, we have
\begin{align*}
 \int&\bbone_{\clP(\Span(t_B))}(t_{[l-1]}) \, \dd\!\sigma_l(t) 
=  \int\bbone_{\clP(\Span(t_B))}(t_{[l-1]}) \, \dd\!\sigma_{l,l}(t_c \mid c \in \Delta_0 S)\\
&=\int \left( \int \bbone_{\clP(\Span(t_B))}(t_{[l-1]}) \, \dd\!\lambda_\perp^{(t_c \mid c < [l-1])}(t_{[l-1]}) \right)  \dd\!\sigma_{l,l-1}(t_c \mid c < [l-1]),		
\end{align*}
hence our assertion follows after showing that
\[
	I:= \clP\big(\!\Span(t_B)\big) \cap \clP\big(\!\Span(t_c \mid c < [l-1])^\perp\big) = \emptyset
\]
for almost every $(t_c \mid c < [l-1])$. Let $\Sigma := \Span(t_c \mid c \in \St(B))$, where $\St(B)$ is the subset of vertices of $\Delta(l)$ introduced in \autoref{def:star}. For $i \in [l-1]$, we define
\[
p_i := t_{B_i} \qand q_i:= t_{\mathrm{op}_i(B)}, 
\]
where $\mathrm{op}_i(B)$ is defined as in \autoref{def:star}. If $L := \bbP^{-1}(\Sigma)$, we have
\[
\bbP(L) = \Sigma = \Span(p_0, \dots, p_{l-1}, q_0, \dots, q_{l-1}).
\]
(see \autoref{StarBExplicit}). We have $\clP\big(\!\Span(t_B)\big)  \subset \clP(\Sigma)$ and $\clP(\Sigma) \subset \clP\big(\!\Span(t_c \mid c < [l-1])\big)$, hence $ \clP\big(\!\Span(t_c \mid c < [l-1])^\perp\big) \subset \clP(\Sigma^\perp)$. It thus remains to show only, that $L$ is non-degenerate. For this we verify the conditions of \autoref{linalg}.

 Since $B$ is a simplex in $\Delta(l)$, its image $t_B$ is a simplex in $\BX$ and hence $\omega(p_i, p_j) = \omega(t_{B_i}, t_{B_j}) = 0$ for all $i,j \in [l-1]$ with $i \neq j$. On the other hand, if $i \neq j$, then ${B_i}$ and $\mathrm{op}_j B$ are contained in a common simplex (i.e. $C'$ with $\delta_jC = \delta_m C'$) and hence similarly $\omega(p_i, q_j) = \omega(t_{B_i}, t_{\mathrm{op}_jB}) = 0$. 
 
It remains to show that $p_i =  t_{B_i}$ and $q_i =  t_{\mathrm{op}_i(B)}$ satisfy $\omega(p_i, q_i) \neq 0$. Let $m$ be the smallest integer such that $\{B_i, \mathrm{op}_i(B)\} \subset \Delta_0(l)_{m}$ holds. One can prove by induction that there exists a subset $\Theta_{B,i} \subset \Delta_0(l)_{m}$ and a measure $\xi^{L_0,L_1} \in \Prob(\clP^{\Theta_{B,i}})$ for almost every pair $(L_0,L_1) \in (\clG_{r-1})^2$ such that for every $\varphi \in \scL^1(\clP^2)$, the equality
\begin{align} 
	\label{aaaaa} \int \varphi(t_{B_i},t_{\mathrm{op}_i(B)}) \, \dd\!\sigma_l(t) = \int&  \int \left(\int \int \varphi(t_{B_i},\, t_{\mathrm{op}_i(B)}) \, \dd\!\lambda_{A_1}(t_{B_i}) \, \dd\!\lambda_{A_2}(t_{\mathrm{op}_i(B)}) \right) \\
	& \dd\!\xi^{L_0,L_1}(t_{c'} \mid c' \in \Theta_{B,i}) \, \dd\!\nu^{\otimes 2}(L_0,L_1) \nonumber
\end{align}
holds, with $A_1$ [resp. $A_2$] being either one of the Lagrangians $L_0$, $L_1$, or of the form $\Span(t_c \mid c < B_i)^\perp$ [resp. $\Span(t_c \mid c < \mathrm{op}_i(B))^\perp$]. Now, we let $\varphi$ be the characteristic function of the set $\{(a_1,a_2) \in \clP^2 \mid \omega(a_1,a_2) \neq 0\}$. We claim then that the term in parentheses in the equality above equals 1 almost surely, and this finishes the proof by Fubini. 

To prove the claim, we verify the hypotheses of \autoref{lemmaaaaa}. Recall that $r \geq 5$ (see \autoref{rank_geq_5}) and that the codimensions of $\Span(t_c \mid c < B_i)^\perp$ and $\Span(t_c \mid c < \mathrm{op}_i(B))^\perp$ are bounded above by $r-2$, as observed in \eqref{annoying_bound}. Moreover, since the vertices $B_i$ and $\mathrm{op}_i(B)$ are not in a common simplex in $\Delta(l)$, we have \vspace{-5pt}
\begin{itemize}[leftmargin=20pt]
	\item that $A_1 \neq A_2$ almost surely if both $A_1,A_2$ are Lagrangians, and
	\item that the random variables $t_{B_i}$ and $t_{\mathrm{op}_i(B)}$ are almost surely independent, i.e. the term in parenthesis in \eqref{aaaaa} is an integral with respect to the product measure $\lambda_{A_1} \otimes \lambda_{A_2}$. \qedhere
\end{itemize}
\end{proof}

\begin{proof}[Proof of \autoref{PBSMeasures}] It remains to show only that the measures $(\sigma_l)_{l \in [\gamma+1]}$ constructed in \autoref{SigmaFinal} satisfy Properties (iii), (iv) and (v) of \autoref{PBSMeasures}.

To prove (iv), observe that the diagram below commutes, where the right vertical arrow is the projection $\delta_\emptyset: \BX_l \times \clP \to \BX_l$ onto the first factor.\vspace{-11pt}
\begin{equation*} 
		\begin{gathered}	
		\xymatrixcolsep{3pc}
		\xymatrix@R=14pt{\scD_{l+1} \ar@{<->}[d]_{\cong} \ar[rr]^-{\partial_\emptyset} & & \BX_{l} \\
		\llap{$\clP^{\Delta_0(l+1)} = \,  $}\clP^{\Delta_0(l+1)_{l+1}}\ar[r]^-{\pi_0 \circ \cdots \circ \pi_l} &\clP^{\Delta_0(l+1)_{0}} \ar@{<->}[r]^{\cong} & \BX_{l} \times \clP \ar[u]_{\delta_\emptyset}} 
		\end{gathered}
\end{equation*}
The left-bottom arrow is p.m.p.\ by  \autoref{brute}, and both vertical arrows and the right-bottom one are p.m.p.\ by definition. This implies that the top horizontal arrow is p.m.p., proving (iv).

Similarly for (v), we consider the following commutative diagram, where $\eta_i := \eta_{[l]\smallsetminus \{i\} \, \subset \, [l]}$ and $\varepsilon_{i}$ is the relabelling map from \autoref{ProjectionsEpsiloni}:
\begin{equation*} 
		\begin{gathered}	
		\xymatrixcolsep{3pc}
		\xymatrix@R=18pt{\scD_{l+1} \ar@{<->}[d]_{\cong} \ar[rr]^-{\partial_i} & & \scD_l \\
	\llap{$\clP^{\Delta_0(l+1)} = \,$} \clP^{\Delta_0(l+1)_{l+1}}  \ar[r]^-{\eta_i \circ \pi_{l}} & \clP^{\Delta_0{[l] \smallsetminus \{i\}}_l} \ar[r]_-{\cong}^-{\varepsilon_{i}} & \clP^{\Delta_0(l)_l} \rlap{$\, = \clP^{\Delta_0(l)}$.} \ar@{<->}[u]_{\cong}} 
		\end{gathered} \vspace{-5pt}
\end{equation*}
Again, the left-bottom arrow is p.m.p by \autoref{brute} and all the canonical bijections above are evidently p.m.p., hence (v) follows.

For (iii), we recall that for $l \in [\gamma+1]$ the canonical measure class on $\BX_l$ is the unique $G$-invariant measure class which gives total mass to the subset $X_l \subset \BX_l$. Since
the pushforward $\ev_{C \, \ast}\sigma_l$ is $G$-quasi-invariant, it thus suffices to show that the equality
\begin{equation} \label{eq:X_l}
\sigma_l\left(\left\{t \in \scD_l \ \mid \  t_C \in X_l \right\}\right) = 1 
\end{equation}
holds for all $C \in \Delta_l(l)$. For the proof of \eqref{eq:X_l}, we proceed by induction on $l$. The case $l=0$ is clear since $\BX_0=X_0$. Now, assume that the claim holds for $l - 1 \in [\gamma]$, and let $C \in \Delta_{l}(l)$. Let us write $C := (B < [l-1])$ for $B \in \Delta_{l-1}([l-1] \smallsetminus \{i\})$ with $i \in [l-1]$. By definition of $\sigma_l$, the left-hand side of \eqref{eq:X_l} equals
\[
I:=\int \left(\int\bbone_{\clP \smallsetminus \clP(\Span(t_{B}))}(t_{[l-1]}) \, \dd\!\lambda_\perp^{(t_{c'} \mid c' < [l-1])}(t_{[l-1]})\right) \cdot \bbone_{X_{l-1}}(t_{B})\, \dd\!\sigma_{l,l-1}(t_c \mid c \in \Lambda_{l,l-1}).
\]
The inner integral equals one by \autoref{interseccion}. Thus,
\[
	I = \int \bbone_{X_{l-1}}(t_{(\varepsilon_i)^{-1}(B)})\, \dd\!\sigma_{l-1}(t) = \ev_{(\varepsilon_i)^{-1}(B) \ \ast} \sigma_{l-1}(X_{l-1}),
\]
where $(\varepsilon_i)^{-1}: \Delta_0([l] \smallsetminus \{i\})_l \to \Delta_0(l)_l$ is the inverse of the relabelling $\varepsilon_i$ from \autoref{ProjectionsEpsiloni}. The induction hypothesis then finishes the proof.
\end{proof}

\begin{rem} \label{problems}
The existence of the measures $\sigma_l \in \Prob(\scD_l)$ was only guaranteed up to degree $\gamma + 1$, and this was the only limitation for extending our partial homotopy to higher degrees. In turn, the definition of $\gamma$ was crucially used in the proofs of \autoref{thm_welldefined} and \autoref{interseccion}.
\end{rem}

At this point we have established \autoref{AcyclicityBound}.

\section{The Quillen argument}

\subsection{Proof of \autoref{ThmIntro}}
The purpose of this section is to derive \autoref{ThmIntro} and \autoref{ThmIntroSG} from \autoref{AcyclicityBound} via the measurable Quillen criterion from \cite{DMH0}. 

Let $(\sigma,\varepsilon) \in \{(\id_\kay,-1), \ (\id_\kay,+1), \ (\,\bar\cdot\,,+1) \}$, and fix an admissible parameter $d \in \bbN$ in the sense of \autoref{restrictions_d}. Then we have the chain of block embeddings of standard $(\sigma,\varepsilon)$-formed spaces
\[
	V_{\sigma,\varepsilon}^{0,d} \hookrightarrow V_{\sigma,\varepsilon}^{1,d} \hookrightarrow V_{\sigma,\varepsilon}^{2,d} \hookrightarrow \cdots \hookrightarrow V_{\sigma,\varepsilon}^{r,d} \hookrightarrow \cdots
\]
indexed by the rank $r$. If we write $G_r := \Aut(V_{\sigma,\varepsilon}^{r,d})$, then we also also the corresponding block embeddings of automorphism groups \vspace{-2pt}
\begin{equation*}
	G_0 \overset{\iota_0}{\hookrightarrow} G_1 \overset{\iota_1}{\hookrightarrow} G_2  \overset{\iota_2}{\hookrightarrow} \cdots \overset{\iota_{r-1}}{\hookrightarrow} G_r \overset{\iota_r}{\hookrightarrow} \cdots
\end{equation*}
Depending on the choice of $(\sigma,\varepsilon)$, we retrieve exactly the classical families from \autoref{ThmIntro}: \vspace{3pt}
\begin{equation*} \begin{array}{lrl}
	\kay \in \{\bbC,\bbR\}, (\sigma,\varepsilon) = (\id,-1): \quad & \{1\} < & \Sp_2(\kay) < \cdots < \Sp_{2r}(\kay) < \cdots \\[2pt]
	
	\kay = \bbC, \ (\sigma,\varepsilon) = (\id,+1), \ d = 0: \quad  &\{1\} < & \OO_2(\bbC) < \cdots < \OO_{2r}(\bbC) < \cdots  
	\\[2pt]
	
	\kay = \bbC, \ (\sigma,\varepsilon) = (\id,+1), \ d = 1: \quad &\{\pm1\} < & \OO_3(\bbC) < \cdots < \OO_{2r+1}(\bbC) < \cdots  
	\\[2pt]
	
	\kay = \bbR, \ (\sigma,\varepsilon) = (\id,+1): \quad &\OO(d) < & \OO(d+1,1) < \cdots < \OO(d+r,r)< \cdots \\[2pt]
	
	\kay = \bbC, \ (\sigma,\varepsilon) = (\, \bar \cdot\,,+1): \quad &\UU(d) < & \UU(d+1,1) < \cdots < \UU(d+r,r)< \cdots
\end{array}\end{equation*}

Let us denote by $X(r)$ the Stiefel complex of $V_{\sigma,\varepsilon}^{r,d}$. By \autoref{thm:transitivity} and \autoref{AcyclicityBound},
\begin{itemize}
\item $X(r)$ is boundedly $\gamma(r)$-acyclic with $\gamma(r)$ as in \autoref{AcyclicityBound};
\item $G_r$ acts $(r-1)$-transitively on $X(r)$.
\end{itemize}
To be able to apply our measurable Quillen criterion \cite[Theorem 4.6]{DMH0} we still need to show that the actions are compatible in a suitable sense.

For $r \geq 0$, we let $\{e_{r}, \dots, e_1, h_1,\ldots,h_d, f_1, \dots, f_{r}\}$ be an adapted basis of $V_{\sigma,\varepsilon}^{r,d}$, and define simplices 
\[	
	o_{r,q} := ([e_{r}], \dots, [e_{r-q}]) \in X(r)_q
\]
for every $q \in [r-1]$. Observe that $o_{r,q}$ is a face of $o_{r,q+1}$ for every $q \in [r-2]$. Thus, if $H_{r,q}$ denotes the stabilizer of $o_{r,q}$, we have inclusions
$
	G_r =: H_{r,-1} > H_{r,0} > \cdots > H_{r,\tau(r)}.
$
We shall denote by $\iota_{r,q+1}: H_{r,q+1} \hookrightarrow H_{r,q}$ the corresponding inclusion. A computation shows that each stabilizer $H_{r,q}$ is represented in the adapted basis as
\begin{align} \label{stabilizer_Hrq}
H_{r, q} &= \left\{\left. M= \begin{pmatrix}
			D & \ast & \ast \\[-2pt]   
			   & A      & \ast \\[-2pt]
			   &         & \Theta(D)
			\end{pmatrix} \ \  \right| \begin{array}{l}
								D \in \GL_{q+1}(\kay) \ \mbox{diagonal} \\
								A \in G_{r-q-1}
			                                     \end{array} \right\}< G_{r},\end{align}
where $\Theta: \GL_{k}(\kay) \to \GL_{k}(\kay)$ is the involution defined as $\Theta(A) := Q_{k} \, \sigma(A)^{-1} \, Q_{k}$, and the asterisks correspond to entries conditioned so that the matrix is in $G_r$. The expression of the stabilizers $H_{r,q}$ given in \eqref{stabilizer_Hrq} implies the following lemma. 

\begin{lem} \label{thm:compatibility}
For every $r \geq 0$ and every $q \in \{-1,0, \dots, r-1\}$, there exist an epimorphism
$
	\pi_{r,q}: H_{r,q} \twoheadrightarrow G_{r-q-1} 
$
with amenable kernel and a continuous homomorphic section $\sigma_{r,q}$ of $\pi_{r,q}$ such that for all $r \geq 0$ and $q \in \{-1,0, \dots r-2\}$ the following diagram commutes: \vspace{-3pt}
\begin{equation*}
		\begin{gathered}	
		\xymatrixcolsep{3.5pc}
		\xymatrix@R=12pt{H_{r,q+1} \ar@{^{(}->}[r]^{\iota_{r,q+1}} & H_{r,q} \ar@{->>}[d]^{\pi_{r,q}} \\
		G_{r-q-2} \ar[u]^{\sigma_{r,q+1}} \ar@{^{(}->}[r]^{\iota_{r-q-2}} & G_{r-q-1}} 
		\end{gathered}
\end{equation*}
\end{lem}

\begin{proof}
Consider the homomorphisms $\pi_{r,q}:H_{r,q} \to G_{r-q-1}$ defined as $\pi_{r,q}(M) :=A$, where $M$ and $A$ are as in \eqref{stabilizer_Hrq}. The kernel of $\pi_{r,q}$ is solvable, hence amenable, and the composition
$
\sigma_{r,q} := \iota_{r-1} \, \circ \, \cdots \, \circ \, \iota_{r-q-2}, 
$ 
whose image is contained in $H_{r,q}$, is a continuous section of $\pi_{r,q}$. That the diagram above commutes is a simple verification. 
\end{proof}
In the language of \cite{DMH0} this establishes that $(G_r, X(r))_{r \geq 0}$ is a ``measurable $(\gamma(r), r-1)$-Quillen family''; then \autoref{Quillen} follows from \cite[Theorem 4.6]{DMH0} (with intial parameter $q_0 := 1$ and $\tau(r) := r-1$), and \autoref{ThmIntro} is an immediate consequence.
\subsection{Proof of \autoref{ThmIntroSG}}
We now deduce \autoref{ThmIntroSG} from \autoref{ThmIntro}. For this we fix a type $(\sigma,\varepsilon) \neq (\id_\kay,-1)$, and an admissible parameter $d$. We then set $G_r := \Aut(V_{\sigma,\varepsilon}^{r,d})$ and
\[
	\mathrm{S}G_r := \{g \in G_r \mid \det g = 1\} \lhd G_r.
\]
Then the sequences $\mathrm{S}G_1 < \mathrm{S}G_2 < \dots$ are precisely the sequences appearing in \autoref{ThmIntroSG}. The corollary thus follows from 
\autoref{ThmIntro} and the following two lemmas.
\begin{lem} \label{isom_SG_G}
Let either $\kay = \bbC$ and $(\sigma,\varepsilon) = (\, \bar\cdot,\,+1)$ or $(\sigma,\varepsilon) = (\id_\kay,+1)$ and $d$ be odd. Then
for every $r \in \bbN$ and $q \in \bbN$ the restriction 
$		\mathrm{res}: \Hb^q(G_r) \to \Hb^q(\mathrm{S}G_r)$ is an isomorphism.
\end{lem}
\begin{proof} We fix $r \geq 1$ and abbreviate $G := G_r$ and $\mathrm{S}G := \mathrm{S}G_r$. Set $\bbK := \bbC$ if $\kay = \bbC$ and $(\sigma,\varepsilon) = (\,\bar\cdot,\,+1)$, and $\bbK := \bbR$ otherwise. Let $n := 2r + d$ be the $\kay$-dimension of the formed space $V_{\sigma,\varepsilon}^{r,d}$. We have the diagram \vspace{-5pt}
\begin{equation*}
	\xymatrixrowsep{1.3pc}
	\xymatrix{1 \ar[r] & \mu_n \ar[r] \ar@{^{(}->}[d] & \bbS_\bbK^1 \ar[r]^-{(-)^n} \ar@{^{(}->}[d] & (\bbS_\bbK^1)^n \ar[r] \ar@{^{(}->}[d] & 1\\
	1 \ar[r] & \mathrm{S}G \ar[r] & G \ar[r]^-{\det} & \bbS_\bbK^1 \ar[r] & 1}
\end{equation*}
where $Z(G) = \bbS_\bbK^1$ is the unit sphere in $\bbK$ with respect to the absolute value, $Z(\mathrm{S}G) = \mu_n < \bbS_\bbK^1$ denotes $n$-th roots of unity, and $(\bbS^1_\bbK)^n < \bbS^1_\bbK$ is the subgroup of $n$-th powers. The leftmost and middle vertical arrows are inclusions as scalar matrices; rows are exact and squares commute. Under our assumptions, we have the isomorphism $\mathrm{PS}G := SG/Z(\mathrm{S}G) \cong G/Z(G) =: \mathrm{P}G$ 
after taking quotients, since $(\bbS^1_\bbK)^n = \bbS^1_\bbK$. Thus, for any $q \geq 0$, the diagram
\[
	\xymatrixcolsep{4pc}
		\xymatrix@R=15pt{\Hb^q(\mathrm{S}G) & \ar[l] \Hb^q(G)  \\
		\Hb^q(\mathrm{PS}G) \ar[u]_{\cong} \ar@{=}[r]& \Hb^q(\mathrm{P}G) \ar[u]^{\cong}} 	
\]
in bounded cohomology commutes. The vertical maps are isomorphisms since the projections that induced them have abelian kernels \cite[Cor. 8.5.2]{Monod-Book}. 
\end{proof}

\begin{lem} \label{isom_SG_G2}
For all $r \geq 1$, $q<2r$ and $d \in \bbN$ the restriction 
\[	\mathrm{res}: \Hb^q(\SO(d+r,r)) \to \Hb^q(\SO_0(d+r,r))\]
is an isomorphism.
\end{lem}
\begin{proof}
Our argument is based on the proof of \cite[Cor. 1.4]{Monod-sot}. We write
\begin{equation*}
	\mathrm{S}G^0 := \SO_0(d+r,r) \lhd \SO(d+r,r) =: \mathrm{S}G  
\end{equation*}
Note that $|\mathrm{S}G:\mathrm{S}G^0| = 2$. By the Eckmann--Shapiro lemma for $\Linfty$-modules (see \cite[Prop. 10.1.3]{Monod-Book}), the restriction $\Hb^q(\mathrm{S}G) \to \Hb^q(\mathrm{S}G^0)$ is an isomorphism for all $q < 2r$ if and only if so is the map 
\[
 	\Hb^q(\mathrm{S}G) \to \Hb^q(\mathrm{S}G;\ell^\infty(\{\pm 1\}))
\]
induced by the inclusion of constants $\bbR \hookrightarrow \ell^\infty(\{\pm 1\})$, where the $\mathrm{S}G$-action on $\{\pm 1\} \cong \mathrm{S}G/\mathrm{S}G^0$ is by left-multiplication. After considering the long exact sequence in bounded cohomology \cite[Prop. 8.2.1]{Monod-Book}, the claim follows if
\[
	\Hb^q(\mathrm{S}G;\ell^\infty(\{\pm 1\})/\bbR) = 0 \quad \mbox{ for all } q<2r.
\]
This is a consequence of \cite[Thm. 1.2]{Monod-sot} since $\mathrm{S}G = \SO(d+r,r)$ is connected as an algebraic $\bbR$-group, and $\ell^\infty(\{\pm 1\})/\bbR$, isomorphic to the space of odd functions $\{\pm 1\} \to \bbR$ as a $G$-module, has no non-zero invariant vectors.
\end{proof}
At this point we have established \autoref{ThmIntroSG}.
\begin{rem} From \autoref{isom_SG_G} we obtain that the sequences
\[
\begin{array}{ll}
	\{1\} < \SO_3(\bbC) < \SO_5(\bbC) < \cdots < \SO_{2r+1}(\bbC) < \cdots \\[2pt]
	\SO(d) < \SO(d+1,1) < \SO(d+2,2) < \cdots < \SO(d+r,r)< \cdots  & \mbox{for \emph{odd} } d \in \bbN \\[2pt]
	\SU(d) < \SU(d+1,1) < \SU(d+2,2) < \cdots < \SU(d+r,r)< \cdots  &\mbox{for any } d \in \bbN
\end{array}
\]
have the same stability range as the corresponding families of the corresponding families of general orthogonal or general unitary groups. The stability range of the family
\[\begin{array}{ll}
	\SO(d) < \SO_0(d+1,1) < \SO_0(d+2,2) < \cdots < \SO_0(d+r,r)< \cdots  & \mbox{for \emph{odd} } d \in \bbN
\end{array}\]
is given by
\[
r(q) = \min\{r_0(q), (q+1)/2\},
\]
where $r_0(q)$ is the stability range of the family $\OO(d) < \OO(d+1, 1) < \dots$.
\end{rem}

\end{document}